\documentclass[12pt]{article}
\usepackage[utf8]{inputenc}
\usepackage[english]{babel}
\usepackage{color}
\usepackage{graphicx}
\graphicspath{{./}}
\usepackage{vmumath}
\usepackage{amsmath}
\usepackage{amsthm}
\usepackage{avo_book}
\usepackage{hyperref}
\usepackage{caption}
\usepackage{booktabs}
\usepackage{subcaption}
\usepackage[a4paper,left=15mm,right=15mm,top=30mm,bottom=20mm]{geometry}
\noaftermath
\begin{document}

\bibliographystyle{unsrt}

\title{Enumeration of Unsensed Orientable and Non-Orientable Maps}

\author{
Evgeniy Krasko \qquad  Alexander Omelchenko\\
\small National Research University Higher School of Economics\\
\small Soyuza Pechatnikov, 16, St. Petersburg, 190008, Russia\\
\small\tt \{krasko.evgeniy, avo.travel\}@gmail.com
}

\begin{abstract}
The paper is devoted to enumeration of unsensed maps, unlabelled maps on surfaces up to all their homeomorphisms. We modify the technique developed Mednykh and Nedela to make it suitable for enumerating unsensed orientable and non-orientable maps on a surface of a given genus $g$. We obtain general formulas which reduce the problem of counting unsensed maps to the problem of enumerating rooted quotient maps on orbifolds. We also derive recurrence relations for quotient maps on the disk and on the M\"obius band, two typical examples of such orbifolds with boundary. In the final part of the paper we describe an algorithm that allowed us to calculate the numbers of unsensed maps on orientable surfaces of genus up to $10$ and on non-orientable surfaces of genus up to $13$ by the number of edges.
\end{abstract}

\maketitle

\section*{Introduction}

By a topological map $\M$ on a surface $X$ we will mean a $2$-cell imbedding of a connected graph $G$, loops and multiple edges allowed, into a compact connected $2$-dimensional manifold $X$ without boundary, such that the connected components of $X-G$ are $2$-cells. The $0$-, $1$-, and $2$-dimensional cells of a map $\M$ are its vertices, edges, and faces, respectively \cite{Handbook_of_Graph_Theory}. In this paper we consider both orientable and non-orientable surfaces without boundary. Every such surface can be characterized by its genus $g$. The orientable surface $X$ of genus $g$ is a sphere with $g$ handles. The non-orientable surface of genus $g$ is a sphere with $g$ holes glued with crosscaps (or Mobius bands). Sometimes instead of $g$ we will use the Euler characteristic $\chi$ of the surface $X$, equal to $2-2g$ in the case of an orientable surface $X^+$ and $2-g$ in the case of a non-orientable surface $X^-$.

Two topological maps $\rm M_1$ and $\rm M_2$ on a surface $X$ are said to be isomorphic if there is a homeomorphism $h$ of $X$ that induces an isomorphism of the underlying graphs $G_1$ and $G_2$. Map isomorphism splits the set of all maps on $X$ into equivalence classes, and each such class is called an unlabelled map. For orientable surfaces we have two types of homeomorphisms, orientation-preserving and orientation-reversing. Unlabelled maps on an orientable surface $X^+$ up to only orientation-preserving homeomorphisms are called sensed maps. Unlabelled maps on an orientable or a non-orientable surface up to all homeomorphisms are called unsensed maps. The main aim of our work is to enumerate unsensed maps on both orientable and non-orientable surfaces.

A general technique for counting sensed planar maps on the sphere was developed by Liskovets \cite{Liskovets_85} in the early eighties. His approach reduces the enumerating problem for sensed maps on the sphere to counting {\em quotient maps on orbifolds}, maps on quotients of a surface under a finite group of automorphisms. His ideas were further developed in a series of papers devoted to enumeration of sensed maps on orientable surfaces of a given genus $g$ \cite{Mednykh_Nedela}, regular sensed maps on the torus \cite{Torus_Part_I}, regular sensed maps on orientable surfaces of a given genus $g$ \cite{POMI_Reg_maps_english}, sensed hypermaps \cite{Mednykh_Hypermaps}, one-face regular sensed maps \cite{Krasko_Omelch_4_reg_one_face_maps} and one-face maximal unsensed maps \cite{Krasko_POMI}. In all these papers a geometric approach based on enumeration of rooted quotient maps on cyclic orbifolds was employed. A similar approach was also employed in  \cite{Azevedo} to enumerate unsensed orientable maps on surfaces regardless of genus. The authors of \cite{Azevedo} pointed out that their technique can be used to enumerate maps by genus as well. 

In the present paper we modify the technique developed in the papers \cite{Mednykh_Nedela}, \cite{Mednykh_Hypermaps} and \cite{Azevedo} to make it suitable for enumerating unsensed maps on orientable and non-orientable surfaces of a given genus $g$. We obtain two general formulas (\ref{eq:unsensed_orientable_final}) and (\ref{eq:unsensed_nonorientable_final}) in the form of linear combinations of numbers of quotient maps on cyclic orbifolds with integer coefficients. These coefficients are expressed through the numbers of epimorphisms from fundamental groups $\pi_1(O)$ of orbifolds to cyclic groups $\Z_d$, for which we obtain exact analytical expressions. Having these formulas one can enumerate different types of unsensed maps on surfaces of a given genus $g$ (arbitrary maps, regular maps, one-face maps, etc.). In particular, in \cite{Torus_Part_II} we already used this approach to obtain analytical formulas for the numbers of $r$-regular maps on the torus. 

The formulas (\ref{eq:unsensed_orientable_final}) and (\ref{eq:unsensed_nonorientable_final}) reduce the problem of enumerating unsensed maps to the problem of enumerating quotient maps on orbifolds. In the present paper we derive recurrence relations for the numbers of such quotient maps for two typical examples of orbifolds with boundary, the disk and the M\"obius band. Unfortunately, for arbitrary orbifolds such enumerating problem is unlikely to have an analytical solution. It turns out that these orbifolds can be surfaces with $r$ branch points, $h$ boundary components and $g$ handles or cross-caps, and the recurrence relations for the numbers of quotient maps on them can depend on a large number of additional parameters. Nevertheless, we managed to analyze all possible cases and construct an algorithm that allowed us to calculate the total number of unsensed maps on orientable surface of genus up to $10$ and on non-orientable surfaces of genus up to $13$ with a given number of edges. To the best of our knowledge there are no published analytical results on enumerating unsensed orientable maps for such values of $g$ and no corresponding results for non-orientable maps for any $g$.

\section{The basic principles of sensed and unsensed map enumeration}

The Burnside's lemma is typically used as the main tool for enumerating combinatorial objects up to their symmetry group (see, for example, \cite{Chalambides}). This lemma reduces the problem of enumerating such objects to enumeration of labelled objects that have a trivial symmetry group. As it was noted in \cite{Tutte_Census}, for maps on surfaces it is convenient to consider so-called rooted maps as labelled objects. 

A map is called rooted if one of its edges is distinguished, oriented, and assigned a left and a right side (see, for example, \cite{Handbook_of_Graph_Theory}, \cite{Walsh_Lehman}). For enumerating maps on orientable surfaces it is sufficient to distinguish one edge-end, called a dart. In the papers \cite{Tutte_Census}, \cite{Tutte_triangulations}, \cite{Tutte_slicings} W.~Tutte proposed a quite simple and effective method of obtaining recurrence relations for the numbers of rooted maps on a sphere based on the idea of removing or contracting a root edge. In  \cite{Walsh_Lehman} Walsh and Lehman generalized this approach and obtained some recurrence relations for the numbers of rooted maps on arbitrary oriented surfaces of genus $g$. In the end of the eighties Bender and Canfield (see \cite{Bender_asymptotic_number}, \cite{Bender_1988}) obtained analytical solutions of such recurrence relations for rooted maps on some of the simplest surfaces. Some recent results regarding enumeration of rooted maps on surfaces of higher genera can be found in \cite{Walsh_Giorgetti_2014} and \cite{Walsh_Giorgetti}.

Let $X_g^+$ be a closed orientable surface of genus $g$. In the paper \cite{Mednykh_Nedela} with the help of the Burnside's lemma and some additional algebraic and topological considerations Mednykh and Nedela derived the following important formula for determining the numbers $\tilde{\tau}_{X_g^+}(n)$ of sensed orientable maps with $n$ edges, i.e. for counting isomorphism classes of orientable maps up to orientation-preserving homeo\-morphisms: 
\begin{equation}
\label{eq:Mednykh_Ned_2006}
\tilde{\tau}_{X_g^+}(n)=\dfrac{1}{2n}\sum_{\substack{l\mid 2n\\l\cdot m=2n}}\sum_{O\in {\rm Orb}(X_g^+/\Z_{l})}
\Epi_o(\pi_1(O), \Z_{l}) \cdot \tau_O(m).
\end{equation}
Here $O=X_g^+/\Z_l$ is a quotient of the surface $X_g^+$ under the action of a cyclic subgroup $\Z_l$ of the group of automorphisms of $X_g^+$, $\Epi_o(\pi_1(O), \Z_{l})$ are integer coefficients equal to the numbers of order-preserving epimorphisms from the fundamental group $\pi_1(O)$ of the orbifold $O$ onto the cyclic group $\Z_{l}$, and $\tau_O(m)$ are the numbers of rooted quotient maps with $m$ darts on a cyclic orbifold $O$. 

\begin{figure}[ht]
\centering
	\begin{subfigure}[b]{0.4\textwidth}
	\centering
    		\includegraphics[scale=1.3]{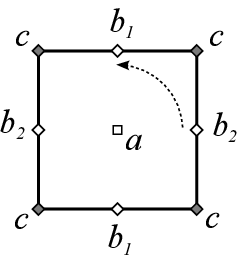}
		\caption{}
	\end{subfigure}	
	\begin{subfigure}[b]{0.4\textwidth}
	\centering
    		\includegraphics[scale=1.3]{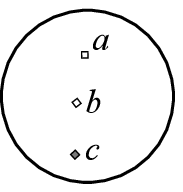}
		\caption{}
	\end{subfigure}	
	\caption{A homeomorphism of the torus and the corresponding orbifold}
\label{fig:orbifold}
\end{figure}

To illustrate these concepts consider the representation of a torus $T=X_1^+$ as a square with its opposite sides identified pairwise (Figure \ref{fig:orbifold} (a)). Rotation of this square by $90^\circ$ ($l=4$) is a typical example of a periodic orientation-preserving homeomorphism of the torus. This homeomorphism splits the set of its points into two subsets, an infinite set of points in the general position and a finite set of singular points (see Figure \ref{fig:orbifold} (a)). Points in the general position are those that lie on some orbit of length $l=4$. Singular points are the remaining ones, and they necessary lie on orbits of smaller length. In our example there are four singular points: $a$, $c$, $b_1$ and $b_2$. The former two of them are fixed, and the latter two are transformed into each other by the rotation by $90^\circ$. Gluing together the points of each orbit of the rotation, we obtain an orbifold $O$, in this case a sphere (Figure \ref{fig:orbifold} (b)). Critical points on the torus get transformed into {\em branch points} of the orbifold $O$ (points $a, b, c$ in Figure \ref{fig:orbifold} (b)). From the topological point of view the described homeomorphism generates a $4$-fold branched covering of the sphere by the torus $T$, and the orbifold $O$ is a quotient $T/\Z_4$.  

In the general case, an orientation-preserving homeomorphism generates the orbifold $O=X_g^+/\Z_l$ which is a surface of genus $\mathfrak{g}$ with a finite number $r$ of branch points. Such orbifold is usually described by its signature
$$
O(\mathfrak{g} ,[m_1,\ldots,m_r]),\qquad 1<m_1\leq\ldots\leq m_r,
$$
where $m_i$ are the {\em branch indices} of the corresponding branch points; each $m_i$ is equal to the period $l$ of the homeomorphism divided by the number of preimages of the corresponding branch point. For the example of the orbifold $O$ shown in Figure \ref{fig:orbifold}, the branch points $a$ and $c$ have branch indices equal to $4$ and the branch index of $b$ is equal to $2$. Consequently, the signature of the corresponding orbifold takes the form
$$
O(0;[2,4,4])\equiv O(0;[2,4^2]).
$$
For the case of the torus there is also one more periodic homeomorphism that preserves its orientation and yields an orbifold with the same signature: the rotation of the square by an angle of $270^\circ$. The coefficient $\Epi_o(\pi_1(O), \Z_{l})$ in (\ref{eq:Mednykh_Ned_2006}) is responsible for calculating all homeomorphisms leading to the same orbifold $O$.

\begin{figure}[ht]
\centering
	\begin{subfigure}[b]{0.4\textwidth}
	\centering
    		\includegraphics[scale=1.3]{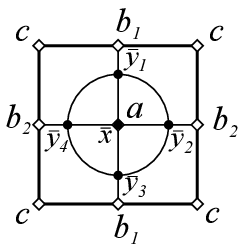}
		\caption{}
	\end{subfigure}
	\begin{subfigure}[b]{0.4\textwidth}
	\centering
    		\includegraphics[scale=1.3]{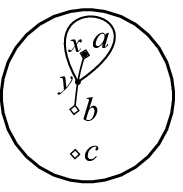}
		\caption{}
	\end{subfigure}
\caption{}
\label{fig:orbifold_map}
\end{figure}

To illustrate the concept of a quotient map on an orbifold $O$, consider, for example, a map $\M$ on the torus $T$ which is symmetric under the rotation of the square by $90^\circ$ (Figure \ref{fig:orbifold_map} (a)). Identifying all points lying on each orbit of the rotation we obtain the quotient map $\mathfrak{M}$ on $O$, shown in the Figure \ref{fig:orbifold_map} (b). This quotient map would be a map on the sphere with the numbers of vertices, edges and faces equal to those of the original map $\M$ divided by $4$ if the orbifold $O$ had no branch points and the surface $X_1^+=T$ had no corresponding critical points. The existence of these points makes this correspondence more complicated.

Assume that a vertex $x$ of a quotient map $\mathfrak{M}$ coincides with some branch point of an index $m_i$ on the orbifold $O$ (see vertex $x$ in Figure \ref{fig:orbifold}(b) which coincides with $a$). Then this vertex corresponds to $l/m_i$ vertices of the map $\M$ on the original surface $X_g^+$. The degree $d$ of $x$ in this case gets multiplied by $m_i$ on $X_g^+$ and becomes equal to $m_id$. For example, the vertex $x$ of degree $1$ of the quotient map $\mathfrak{M}$ shown in Figure \ref{fig:orbifold} (b) corresponds to a single vertex $\bar{x}$ of degree $4$ of the map $\M$ in Figure \ref{fig:orbifold}(a).

Now assume that a branch point of an orbifold $O$ falls into some face $f$ of the quotient map $\mathfrak{M}$ (see branch point $c$ in Figure \ref{fig:orbifold_map}(b)). This point will correspond to $l/m_i$ points on the surface $S_g$, $m_i$ being the corresponding branch index. The remaining points of $f$ are not branch points, so each of them corresponds to $l$ points on $S_g$. Hence, as in the case of a vertex, the degree of the face $f$ is multiplied by $m_i$ when this face is lifted to the surface $S_g$. For example, the degree of the face that contains the branch point $c$ in Figure \ref{fig:orbifold_map}(b) is multiplied by $4$ on the torus $T$.

One more property of quotient maps is the possibility of having {\em dangling semiedges} which end not in vertices but in branch points of degree $2$ (see branch point $b$ in Figure \ref{fig:orbifold_map}(b)). When lifted to the surface $X_g^+$, any such edge gets transformed into $l/2$ edges of $\M$. These edges on the surface $X_g^+$ contain critical points $b_i$ corresponding to the branch point $b$ (see critical points $b_1, b_2$ in Figure \ref{fig:orbifold_map} (a)). If an orbifold $O$ has no branch points of index $2$, then there are no dangling semiedges in any quotient map $\mathfrak{M}$ on $O$. 

Summing up, we can say that to use the formula (\ref{eq:Mednykh_Ned_2006}) we have to solve three subproblems. First, we need to describe the set $X_g^+/\Z_l$ of all orbifolds for a given surface $X_g^+$ and a fixed number $l$. Second, we need to calculate the numbers $\Epi_o(\pi_1(O), \Z_{l})$ of order-preserving epimorphisms from the fundamental group $\pi_1(O)$ of the orbifold $O$ onto the cyclic group $\Z_{l}$. Finally, we need to enumerate all rooted quotient maps on each orbifold $O$. The solutions of all these three subproblems were described in the paper \cite{Mednykh_Nedela}.

\begin{figure}[ht]
\centering
	\centering
    	\includegraphics[scale=0.7]{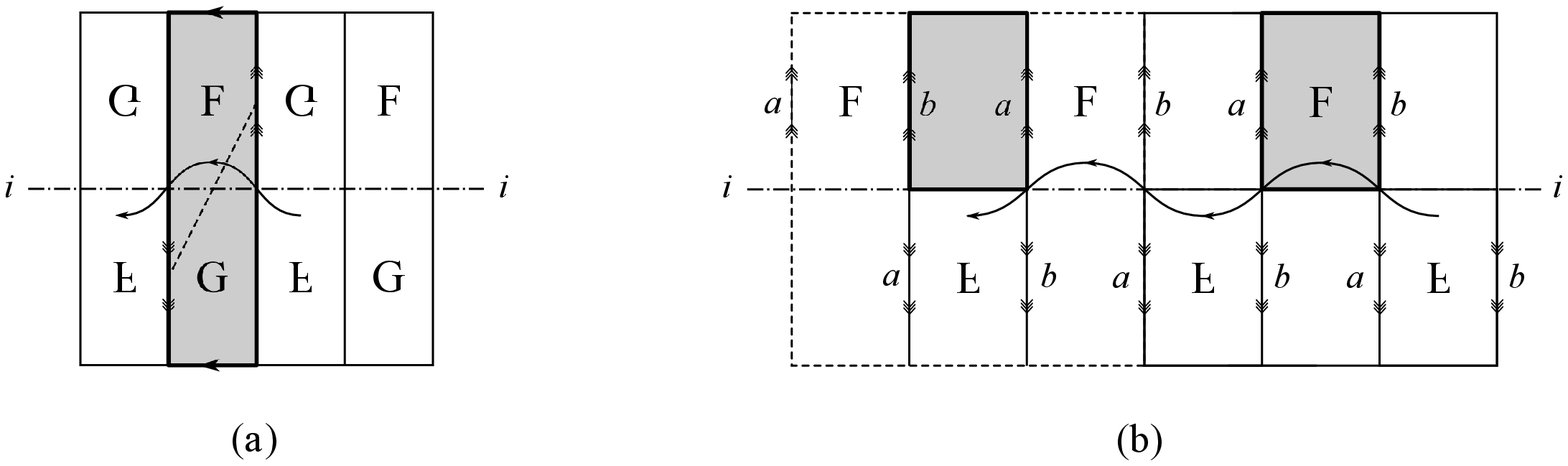}
	\caption{}
\label{fig:Glide_reflection_horis}
\end{figure}

Next we proceed with unsensed maps. Recall that by an unsensed map we mean a equivalence class of isomorphic maps (orientable or not) where the equivalence relation is given by arbitrary map homeomorphisms. In the case of orientable maps, along with the homeomorphisms which preserve the orientation, we must also take orientation-reversing homeomorphisms into account. For a torus the simplest example of such homeomorphism is a glide reflection with respect to the horizontal axis $i$ of the square representing this torus on the plane (see Figure \ref{fig:Glide_reflection_horis}). Let the ratio between the value of the `shift' and the length of the side be a rational number $p/q$, $0\leq p/q<1$, $p$ and $q$ are coprime. In Figure \ref{fig:Glide_reflection_horis} an example of a glide reflection with respect to a horizontal axis and $p/q=1/4$ is shown. The fundamental polygon in this case is one fourth of a square. Since the right side of this polygon is transformed into its left side under the glide reflection with respect to $i$, we may think of these sides as glued together in the reverse direction. Consequently, such homeomorphism generates a $4$-fold branched covering of the Klein bottle $O$ by the torus.

Now consider an example of a glide reflection for $p/q=1/3$ (Figure \ref{fig:Glide_reflection_horis} (b)). For this ratio of $p$ and $q$ the fundamental polygon is one sixth of the square (shaded area on Figure \ref{fig:Glide_reflection_horis}(b)). Indeed, it would take six steps for the glide reflection to transform each point of the torus into itself. After the second step the left side $a$ of the fundamental polygon will coincide with its right side $b$, and vice versa after the fourth step. At the same time its top and bottom sides will never become coincident. Consequently this glide reflection corresponds to a rectangular fundamental region with its right and left sides glued together. In other words, in this case the orbifold $O$ is an annulus.

\begin{figure}[h]
\centering
\includegraphics[scale=0.7]{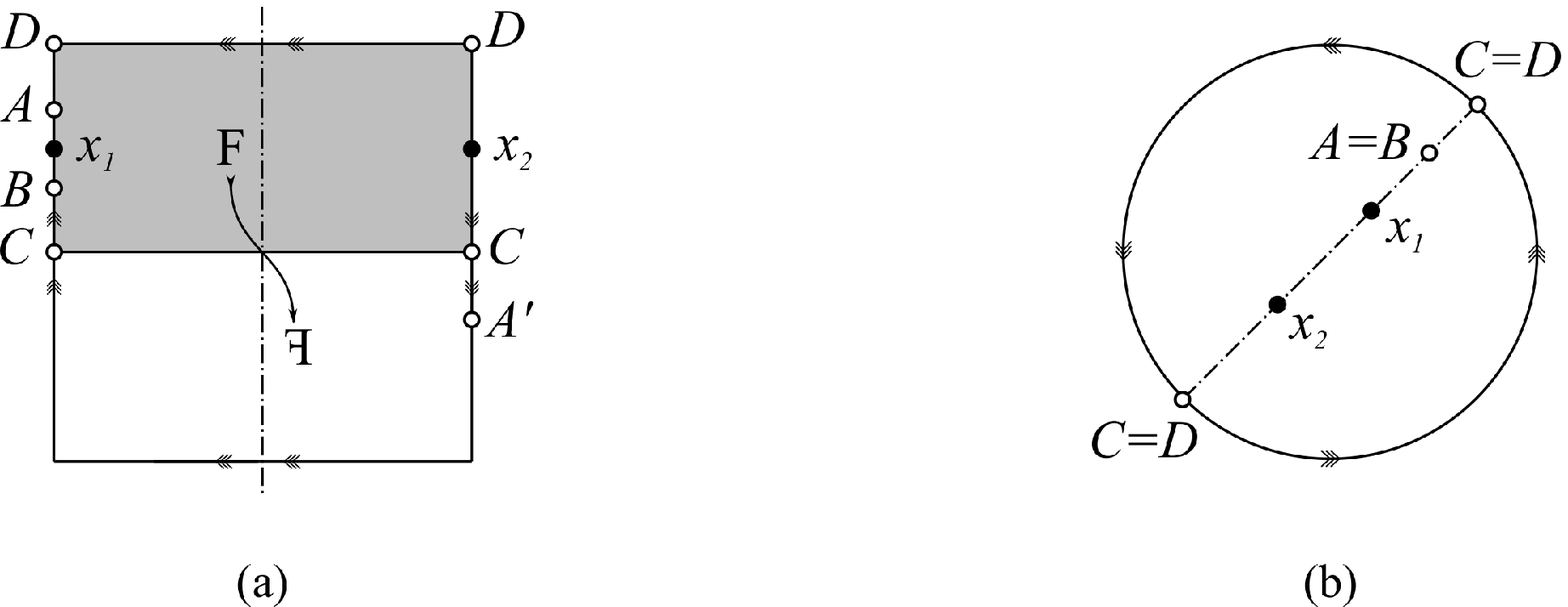}
\caption{}
\label{fig:Proj_plane_2_branch_points}
\end{figure}

As a typical example of a non-orientable surface $X_{g}^-$ consider the Klein bottle. We will use the representation of the Klein bottle as a square with its top and bottom sides glued in the forward direction, and its left and right sides glued in the reverse direction (see Figure \ref{fig:Proj_plane_2_branch_points}(a)). The homeomorphism of this surface that shifts the upper half of the square down with a flip relative to the vertical axis leads to a projective plane with two branch points as the orbifold (Figure \ref{fig:Proj_plane_2_branch_points}(b)). Indeed, consider the cell $F$ shaded in the Figure \ref{fig:Proj_plane_2_branch_points}(a). The point $A$ lying on the upper left boundary of the cell under the action of this homeomorphism is transformed into the point $A'$ on the right border. But this point in the Klein bottle coincides with the point $B$. In other words, for such homeomorphism the points $A$ and $B$ on the left boundary of the cell $F$ are glued together. Points on the right border of the square behave similarly. At the same time, the points $x_1$ and $x_2$ in Figure \ref{fig:Proj_plane_2_branch_points} are transformed into themselves. As a consequence, these points correspond to branch points of index $2$ on the orbifold (Figure \ref{fig:Proj_plane_2_branch_points}(b)). It remains to note that the top and bottm boundaries of the cell $F$ are glued together in the opposite direction. Representing them as a border of a circle (Figure \ref{fig:Proj_plane_2_branch_points}(b)) we obtain a projective plane with two branch points.

\begin{figure}[h]
\centering
\includegraphics[scale=0.7]{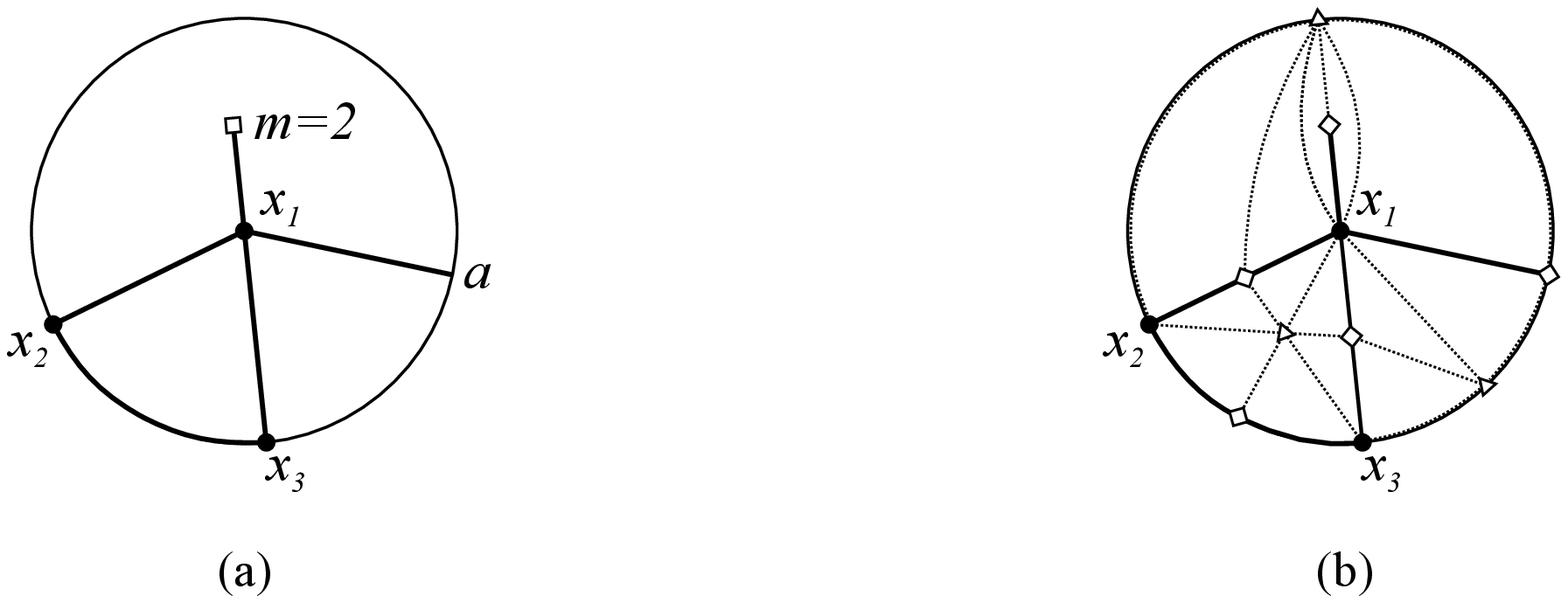}
\caption{}
\label{fig:disk_flags_example}
\end{figure}

As it can be seen from these examples, in the case of arbitrary homeomorphisms the orbifolds can be both orientable and non-orientable surfaces, both closed surfaces and surfaces with boundary. Some branch points may lie on these orbifolds as well, so for arbitrary periodic homeomorphisms we must be able to enumerate quotient maps not only on closed orientable orbifolds, but also on orientable or non-orientable orbifolds with boundary and branch points. Quotient maps on such orbifolds, apart from complete edges and semiedges (see a complete edge $x_1x_2$ in Figure  \ref{fig:disk_flags_example}(a) and a semiedge going from $x_1$ to the branch point of the index $m=2$ in this Figure) may also have so-called halfedges (see halfedge $x_1a$ in Figure \ref{fig:disk_flags_example}(a)) and boundary edges (the edge $x_2x_3$ in Figure \ref{fig:disk_flags_example}(a)). If we place new vertices into the centers of edges (see squares in Figure \ref{fig:disk_flags_example}(b)) and into the centers of faces (see triangles in Figure \ref{fig:disk_flags_example}(b)), and connect neighboring vertices by new edges (see dashed lines in Figure  \ref{fig:disk_flags_example}(b)), we will obtain a partition of the map into triangles which are called flags of this map. Since each edge is incident to an even number of flags regardless of its type (complete edge, semiedge etc.), the total number of flags is even for any quotient map of any orbifold. 

\section{Formulas for counting of unsensed orientable and non-orientable maps}

In the paper \cite{Torus_Part_II} using some direct topological considerations was obtained the following formula for calculating unsensed maps on the torus:
\begin{equation}
\label{eq:r_reg_maps_torus}
\begin{aligned}
\bar{\tau}(n)&=
\dfrac{\tilde{\tau}(n)}{2}+
\dfrac{1}{4n}\sum\limits_{(2d+1)\mid n}\phi(2d+1)\cdot \left[C(2n/(2d+1))+M(2n/(2d+1))\right]+\\
&+\dfrac{1}{4n}\sum\limits_{(2d+1)\mid n}\phi(2d+1)\cdot \kappa(n/(2d+1))+\dfrac{1}{n}\sum\limits_{2d\mid n}\phi(2d)\cdot \kappa(n/(2d)).
\end{aligned}
\end{equation}
Here $\tilde{\tau}(n)$ is the number of sensed maps on the torus, $\kappa(n)$, $C(n)$ and $M(n)$ are the numbers of quotient maps with $n$ darts on the Klein bottle, the annulus and the M\"obius band, respectively. For counting unsensed maps in the general case we must have an analogous formula expressing the number $\bar\tau_{X_g}(n)$ of such maps in terms of numbers of sensed maps and numbers of quotient maps on orbifolds described above.

The same arguments that were used in the proof of formula (\ref{eq:Mednykh_Ned_2006}) in the paper \cite{Mednykh_Hypermaps} (see Theorem 3.3 and its proof in \cite{Mednykh_Hypermaps}) allow us to generalize the formula (\ref{eq:r_reg_maps_torus}) and prove the following statement.

\begin{theor_eng}\label{theor:Med_Ned_generalization}
The number $\bar\tau_{X_\chi}(n)$ of unsensed maps with $n$ edges on surfaces $X_\chi$ of Euler characteristic $\chi$ is calculated by the formula 
\begin{equation}
\label{eq:Med_Ned_generalization}
\bar\tau_{X_\chi}(n) = \frac{1}{4n} \sum\limits_{\substack{m|4n\\lm = 4n}} \sum_{O \in \Orb(X_\chi/\Z_l)} \tau_O(m) \cdot \Epi_o(\pi_1(O), \Z_l).
\end{equation}
Here $\Orb(X_\chi/\Z_l)$ is the set of all orbifolds $O$ arising as quotients of a surface $X_\chi$ by the cyclic group $\Z_l$ of automorphisms of $X_\chi$, $\tau_O(m)$ is the number of rooted quotient maps with $m$ flags on an orbifold $O$,
$\Epi_o(\pi_1(O), \Z_l)$ is the number of order-preserving epimorphisms from the fundamental group $\pi_1(O)$ of the orbifold $O$ onto the cyclic group $\Z_l$.
\end{theor_eng}

In our case the fundamental group $\pi_1(O)$ of the orbifold $O$ is a NEC-group (see \cite{NEC-groups}, \cite{Etayo}, \cite{Non-orient_cyclic_autom}). For orientable orbifolds this group is generated by elliptic elements $x_1,\ldots,x_r$ (loops which go around the branch points), orientation-preserving elements $e_1,\ldots,e_h$ (loops that go around boundary components), reflections $c_1,\ldots,c_h$ and hyperbolic translations $a_1,b_1$, $\ldots$, $a_{\mathfrak{g}},b_{\mathfrak{g}}$, subject to the following relations (see \cite{NEC-groups}, \cite{Etayo}):
\begin{equation}
\label{eq:def_relations_orient}
\prod\limits_{i=1}^{r}x_i\prod\limits_{j=1}^{h}e_j\prod\limits_{k=1}^{\mathfrak{g}}[a_k,b_k]=1,\qquad\qquad 
x_i^{m_i}=1,\qquad\qquad c_j^2=1,\qquad\qquad [c_j,e_j]=1.
\end{equation}

For non-orientable orbifolds the list of generators has the form (see \cite{Non-orient_cyclic_autom})
$$
x_1,\ldots,x_r,c_1,\ldots,c_h,e_1,\ldots,e_h,a_1,\ldots,a_{\mathfrak{g}},
$$
and the relations are
\begin{equation}
\label{eq:def_relations_nonorient}
\prod\limits_{i=1}^{r}x_i\prod\limits_{j=1}^{h}e_j\prod\limits_{k=1}^{{\mathfrak{g}}}a_k^2=1,\qquad\qquad 
x_i^{m_i}=1,\qquad\qquad c_j^2=1,\qquad\qquad [c_j,e_j]=1.
\end{equation}

In the case of an odd $\chi$ there is only one surface of a given Euler characteristic, namely, a non-orientable surface. For such values of $\chi$ the formula (\ref{eq:Med_Ned_generalization}) allows us to calculate the number of unsensed maps on such surface at once. In the case of an even $\chi$ there are two surfaces that correspond to the same value of the Euler characteristic $\chi$, one of which is non-orientable, and the other is orientable. So in the case of an even $\chi$ we need some formulas analogous to (\ref{eq:Med_Ned_generalization}) which would count unsensed maps on orientable and non-orientable surfaces separately.

For enumerating unsensed maps on an orientable surface we introduce the concept of a group with sign structure \cite{Azevedo} — a group $G$ together with a homeomorphism $\sigma$ from $G$ to the group of two elements (''signs'') $(\{+1,-1\},\times)$. For orientable orbifolds and $G=\pi_1(O)$ the following relations must be satisfied (see \cite{Azevedo}):
$$
\sigma(x_i)=\sigma(e_j)=\sigma(a_k)=\sigma(b_k)=+1,\qquad\qquad \sigma(c_j)=-1.
$$
For non-orientable orbifolds the corresponding relations have the following form:
$$
\sigma(x_i)=\sigma(e_j)=+1,\qquad\qquad \sigma(c_j)=\sigma(a_k)=-1.
$$
If instead of all epimorphisms in (\ref{eq:Med_Ned_generalization}) we consider only those that preserve the sign structure, then we get the number $\bar\tau_{X_\chi^+}(4n)$ of unsensed maps on the orientable surface $X_g^+$ of Euler characteristic $\chi$:
\begin{equation}
\label{eq:unsensed_orientable}
\bar\tau_{X_\chi^+}(n) = \frac{1}{4n} \sum\limits_{\substack{m|4n\\lm = 4n}} \sum_{O \in \Orb(X_\chi^+/\Z_l)} \tau_O(m) \cdot \Epi_o^+(\pi_1(O), \Z_l)=\frac{1}{4n} \sum\limits_{\substack{m|4n\\lm = 4n}} \sum_{O \in \Orb(X_\chi/\Z_l)} \tau_O(m) \cdot \Epi_o^+(\pi_1(O), \Z_l).
\end{equation}
Here $\Epi_o^+(\pi_1(O), \Z_l)$ is the number of epimorphisms from the NEC-group $\pi_1(O)$ of an orbifold $O$ to the group $\Z_l$, preserving both the order of elements of finite order and the sign of all elements of the group $\pi_1(O)$. The last transition in formula (\ref{eq:unsensed_orientable}) is correct because in the case of non-orientable surfaces $X_\chi^-$ the number $\Epi_o^+(\pi_1(O), \Z_l)$ of sign-preserving epimorphisms is equal to zero. 

Subtracting the formula (\ref{eq:unsensed_orientable}) from the formula (\ref{eq:Med_Ned_generalization}) we obtain the number of  unsensed maps on the non-orientable surface $X_\chi^-$ of Euler characteristic $\chi$:
$$
\bar\tau_{X^-_\chi}(n) = \frac{1}{4n} \sum\limits_{\substack{m|4n\\lm = 4n}} \sum_{O \in \Orb(X_\chi/\Z_l)} \tau_O(m) \cdot (\Epi_o(\pi_1(O), \Z_l) - \Epi_o^+(\pi_1(O), \Z_l))=
$$
\begin{equation}
\label{eq:unsensed_nonorientable}
= \frac{1}{4n} \sum\limits_{\substack{m|4n\\lm = 4n}} \sum_{O \in \Orb(X^-_\chi/\Z_l)} \tau_O(m) \cdot (\Epi_o(\pi_1(O), \Z_l) - \Epi_o^+(\pi_1(O), \Z_l)).
\end{equation}

Returning to the formula (\ref{eq:unsensed_orientable}), we can simplify it by extracting the term (\ref{eq:Mednykh_Ned_2006}) corresponding to orientable sensed maps. To do this, divide all orbifolds into two subclasses: the subclass $\Orb^+(X_\chi^+/\Z_l)$ of orbifolds arising from orientation-preserving and the subclass $\Orb^-(X_\chi^+/\Z_l)$ of orbifolds arising from orientation-reversing homeomorphisms of $X_\chi^+$:
\begin{equation}
\label{eq:unsensed_orientable_two_summands}
\bar\tau_{X_\chi^+}(n) = \frac{1}{4n} \sum\limits_{\substack{m|4n\\lm = 4n}} 
\Bigl[\sum_{O \in \Orb^+(X_\chi^+/\Z_l)} \tau_O(m) \cdot \Epi_o^+(\pi_1(O), \Z_l) + 
\sum_{O \in \Orb^-(X_\chi^+/\Z_l)} \tau_O(m) \cdot \Epi_o^+(\pi_1(O), \Z_l) \Bigr].
\end{equation}
All orbifolds of the first subclass are orientable surfaces possibly with branch points, but without boundaries \cite{Mednykh_Nedela}. As a consequence, for such orbifolds the generators $c_j$ in (\ref{eq:def_relations_orient}) are vanished $(h=0)$, the sign structure $\pi_1(O)$ is trivial, and therefore $\Epi_o^+(\pi_1(O), \Z_l) = \Epi_o(\pi_1(O), \Z_l)$. For orbifolds of the second subclass the period $l$ of the homeomorphism must be even: an orientation-reversing homeomorphism applied an odd number of times reverses the orientation of the surface and as a result is not identity. 

The last simplification common to all considered quotient maps is the parity of the numbers $m$ of flags. Taking into account all these considerations, we rewrite the formula (\ref{eq:unsensed_orientable_two_summands}) as
$$
\bar\tau_{X_\chi^+}(n) = \frac{1}{2} \biggl[ \frac{1}{2n} \sum\limits_{\substack{m|2n\\lm = 2n}} \sum_{O \in \Orb^+(X_\chi^+/\Z_l)} \tau_O(2m) \cdot \Epi_o(\pi_1(O), \Z_l) + 
$$
$$
+ \frac{1}{2n} \sum\limits_{\substack{m|2n\\lm = n}} \sum_{O \in \Orb^-(X_\chi^+/\Z_{2l})} \tau_O(2m) \cdot \Epi_o^+(\pi_1(O), \Z_{2l}) \biggr].
$$
Note that the first term in brackets is the number $\tilde\tau_{X_g^+}(n)$ of sensed maps with $n$ edges on the orientable surface of genus $g=(2-\chi)/2$ (formula (\ref{eq:Mednykh_Ned_2006})). Thus we have proved the following statement generalizing the formula (\ref{eq:r_reg_maps_torus}).

\begin{theor_eng}\label{theor:unsensed_orient}
The number $\bar\tau_{X_\chi^+}(n)$ of unsensed maps with $n$ edges on the orientable surface $X_\chi^+$ of Euler characteristic $\chi$ is calculated by the formula
\begin{equation}
\label{eq:unsensed_orientable_final}
\bar\tau_{X_\chi^+}(n) = \frac{1}{2} \Bigl( \tilde\tau_{X_g^+}(n) + \frac{1}{2n} \sum\limits_{\substack{m|2n\\lm = n}} \sum_{O \in \Orb^-(X_\chi^+/\Z_{2l})} \tau_O(2m) \cdot \Epi_o^+(\pi_1(O), \Z_{2l}) \Bigr).
\end{equation}
\end{theor_eng}

Finally, taking parity of flags in quotient maps on orbifolds into account, we reformulate the theorem  \ref{theor:Med_Ned_generalization} for the case of unsensed maps on non-orientable surface as follows. 

\begin{theor_eng}\label{theor:unsensed_nonorient}
The number $\bar\tau_{X_\chi^-}(n)$ of unsensed maps on the non-orientable surface $X_\chi^-$ of Euler characteristic $\chi$ is calculated by the formula
\begin{equation}
\label{eq:unsensed_nonorientable_final}
\bar\tau_{X_\chi^-}(n) = \frac{1}{4n} \sum\limits_{\substack{m|2n\\lm = 2n}} \sum_{O \in \Orb(X^-_\chi/\Z_l)} \tau_O(2m) \cdot (\Epi_o(\pi_1(O), \Z_l) - \Epi_o^+(\pi_1(O), \Z_l)).
\end{equation}
\end{theor_eng}

Notice that the formula (\ref{eq:unsensed_nonorientable_final}) is correct for both even and odd values of $\chi$: in the case of odd $\chi$ the number $\Epi_o^+(\pi_1(O), \Z_l)$ is equal to zero.

\section{Orbifolds and the corresponding numbers of coverings}

In order to enumerate unsensed maps by the formulas (\ref{eq:unsensed_orientable_final}) and (\ref{eq:unsensed_nonorientable_final}) we should solve three problems: find the sets $\Orb^-(X_\chi^+/\Z_{2l})$ and $\Orb(X_\chi^-/\Z_l)$ of suitable cyclic orbifolds for a given orientable $(X_\chi^+)$ or non-orientable $(X_\chi^-)$ surface, determine the numbers $\Epi_o(\pi_1(O), \Z_l)$ and $\Epi_o^+(\pi_1(O), \Z_l)$ of order-preserving epimorphisms, and find the numbers of quotient maps $\tau_O(2m)$ with $2m$ flags on cyclic orbifolds. In this section we will solve the first two problems. 

The problem of determining the set of cyclic orbifolds for a given surface can be solved using the Riemann–Hurwitz formula
\begin{equation}
\label{eq:Riemann_Hurvitz}
-\chi=l\left[\alpha\mathfrak{g}-2+h+\sum\limits_{i=1}^r\left(1-\dfrac{1}{m_i}\right)\right],
\end{equation}
that connects the Euler characteristics $\chi$ of the original surface $X_\chi$ with the period $l$ of the homeomor\-phism and the parameters of an orbifold $O$. These parameters are as follows. For an orientable orbifold $\alpha=2$ and $\mathfrak{g}$ is the number of handles. For a non-orientable orbifold $\alpha=1$ and $\mathfrak{g}$ is the number of crosscaps. The parameter $h$ defines the number of boundary components. The numbers $m_i$ are branch indices of branch points. 

The Riemann–Hurwitz formula (\ref{eq:Riemann_Hurvitz}) provides a necessary condition that any combination of orbifold parameters must satisfy. To make the set of such combinations to consider finite, we will need some additional constraints on the period $l$ of homeomorphisms. For orientable surfaces of genus $g\geq 2$ (equivalently, $\chi < 0$) it will be sufficient to use the following upper bound (\cite[Corollary 2]{Etayo}):
$$
l\leq 4g+4,\qquad \text{$g$ even;}\qquad\qquad l\leq 4g-4,\qquad \text{$g$ odd.}
$$
For non-orientable surfaces with $g \geq 3$ the corresponding bounds can be found in \cite{Klein_surf_disser}:
$$
l\leq 2g-2,\qquad \text{$g$ even;}\qquad\qquad l\leq 2g,\qquad \text{$g$ odd.}
$$

The four remaining surfaces are the sphere, the projective plane, the torus and the Klein bottle. For them the parameter $l$ can be arbitrarily large, but for enumerating maps with $4n$ flags it is sufficient to assume that $l \leq 2n$, since there must be at least 2 flags in the quotient map.

Having an upper bound for $l$, we can iterate over all suitable values of $\mathfrak{g}$, $h$, $r$ and $m_i$ and hence obtain a finite list of signatures of possible orbifolds which will contain some excess. Then, for each such signature of an orbifold $O$ we should calculate the numbers $\Epi_o(\pi_1(O), \Z_{l})$ and $\Epi_o^+(\pi_1(O), \Z_{2l})$ of epimorphisms from $\pi_1(O)$ to $\Z_{2l}$ and retain only those orbifolds that have a non-zero number of epimorphisms.

The following theorem based on the inclusion-exclusion principle allows to reduce the problem of determining the numbers $\Epi_o(\pi_1(O), \Z_{l})$ and $\Epi_o^+(\pi_1(O), \Z_{2l})$ to a much easier problem of determining the numbers $\Hom_o(\pi_1(O), \Z_{l})$ and $\Hom_o^+(\pi_1(O), \Z_{2l})$ of homomorphisms from $\pi_1(O)$ to $\Z_{l}$ or $\Z_{2l}$.

\begin{theor_eng}[\cite{Azevedo},\cite{Jones_Oxford}]\label{theor:num_epimorph}
The numbers of order-preserving or orientation-and-order-preserving epi\-mor\-phisms can be expressed through the corresponding numbers of homomorphisms by the formulas
$$
\Epi_o(\pi_1(O), \Z_{l}) = \sum\limits_{d|l} \mu(l/d) \cdot {\Hom}_o(\pi_1(O), \Z_{d}),
$$
$$
\Epi_o^+(\pi_1(O), \Z_{2l}) = \sum\limits_{d|l,\,\, l/d \,\,\text{odd}} \mu(l/d) \cdot {\Hom}_o^+(\pi_1(O), \Z_{2d}).
$$
\end{theor_eng}

We should consider three different cases for an orbifold $O$: the case of an orientable surface with $h>0$ boundary components, the case of a non-orientable surface with $h>0$ boundary components and the case of a non-orientable surface without boundary. We give the derivation of the formulas for the first case in detail. 

\begin{propos_eng}\label{propos:num_hom}
The number of order-preserving homomorphisms in the case of an orientable orbifold $O$ with $h>0$ boundary components is equal to
\begin{equation}
\label{eq:num_ord_pres_hom_orient}
\Hom_o(\pi_1(O),\Z_{d})=d^{2\mathfrak{g}+h-1}\cdot \prod\limits_{i=1}^r\phi(m_i)
\end{equation}
if $d$ is even and $m_i|d$, and $0$ otherwise. 
\end{propos_eng}

\evidpEng Consider an arbitrary order-preserving homomorphism $\psi\colon \pi_1(O)\to \Z_{d}$ from the fundamental group $\pi_1(O)$ of an orientable surface $X_{\mathfrak{g}}^+$ to the cyclic group $\Z_d$ for some integer $d>0$. Under this homomorphism $\psi$ the first two equalities in (\ref{eq:def_relations_orient}) are transformed into relations
\begin{equation}
\label{eq:psi_e_i_psi_f_j}
\sum\limits_{i=1}^r \psi(x_i)+\sum\limits_{j=1}^h \psi(e_j)= 0,
\qquad\qquad m_i\cdot\psi(x_i)=0.
\end{equation}
From the second equation and the condition that the order of $\psi(x_i)$ is equal to $m_i$ it follows that there are $\phi(m_i)$ possible values for $\psi(x_i)$ and that $m_i$ divides $d$. The third condition $c_j^2=1$ is transformed in $\Z_{d}$ into equality 
\begin{equation}
\label{eq:psi_c_i}
2\cdot\psi(c_j)=0.
\end{equation}
Since the homomorphism $\psi$ is order-preserving, the order of $\psi(c_j)$ is equal to $2$. Consequently, $d$ is even and $\psi(c_j)$ can take only one value, namely $d/2$. The fourth condition in (\ref{eq:def_relations_orient}) is satisfied trivially in $\Z_d$.

Consider the first equation in (\ref{eq:psi_e_i_psi_f_j}) again. From this equation it follows that $\psi(e_j)$, $j=1,\ldots,h-1$, can be chosen arbitrarily and $\psi(e_h)$ will be uniquely determined by this choice and the choice of $\psi(x_i)$. So the total number of ways to choose all values of $\psi(e_j)$ is equal to $d^{h-1}$. 

Since there are no restrictions on the choice of $\psi(a_k)$ and $\psi(b_k)$, we can choose them in $d^{2\mathfrak{g}}$ ways. Summarizing, we conclude that the number of order-preserving homomorphisms $\psi$ from $\pi_1(O)$ to $\Z_d$ is given by (\ref{eq:num_ord_pres_hom_orient}). \qed

\begin{propos_eng}
The number of orientation-and-order-preserving homomorphisms in the case of an orientable orbifold $O$ with $h>0$ boundary components is equal to
\begin{equation}
\label{eq:num_ord_pres_hom_orient_plus}
\Hom_o^+(\pi_1(O),\Z_{2d})=d^{2\mathfrak{g}+h-1}\cdot \prod\limits_{i=1}^r\phi(m_i)
\end{equation}
if $d$ is odd and $m_i|d$, and $0$ otherwise. 
\end{propos_eng}

\evidpEng The difference between this case and the previous one is that now we need to take the sign structure into account. Namely, instead of an arbitrary order-preserving homomorphism $\psi\colon \pi_1(O)\to \Z_{d}$ we have to consider an order-preserving homomorphism $\psi\colon \pi_1(O)\to \Z_{2d}$ that also preserves the sign structure $\sigma\colon \pi_1(O)\to (\{+1,-1\},\times)$ given by
\begin{equation}
\label{eq:sign_structure_pi_1_O}
\sigma(x_i)=\sigma(e_j)=\sigma(a_k)=\sigma(b_k)=+1,\qquad \sigma(c_j)=-1.
\end{equation}
It is assumed that the sign structure $\sigma'\colon \Z_{2d}\to (\{+1,-1\},\times)$ is given by
$$
\sigma'(2k)=+1,\qqquad \sigma'(2k+1)=-1.
$$

Consider any such homomorphism $\psi\colon \pi_1(O)\to \Z_{2d}$. Observe that the equations (\ref{eq:psi_e_i_psi_f_j}) and (\ref{eq:psi_c_i}) still hold, but now in $\Z_{2d}$. The conditions (\ref{eq:sign_structure_pi_1_O}) mean that $\psi(x_i)$, $\psi(e_j)$, $\psi(a_k)$, $\psi(b_k)$ must all be even. Consequently, they can only be chosen from the subgroup $\Z_{d}$ of even elements of $\Z_{2d}$. Therefore the number of ways to choose these values is still equal to
$$
\prod\limits_{i=1}^r\phi(m_i)\cdot d^{2\mathfrak{g}+h-1}
$$
and all $m_i$ divide $d$. From the condition $\sigma(c_j)=-1$ it follows that $\psi(c_j)$ is odd. Along with the condition (\ref{eq:psi_c_i}) it means that $d$ must be odd and that all $\psi(c_j)$ are uniquely determined.   \qed

\begin{conseq_eng}
The numbers of order-preserving and orientation-and-order-preserving epi\-mor\-phisms for the case of an orientable orbifold $O$ with $h>0$ boundary components are equal to
\begin{equation}
\label{eq:epi_orient_holes}
{\rm Epi}_o(\pi_1(O),\Z_{l})=(m')^{2\mathfrak{g}+h-1}\cdot J_{2\mathfrak{g}+h-1}\Bigl(\dfrac{l}{m'}\Bigr) \cdot \prod\limits_{i=1}^r\phi(m_i),\quad m'=\lcm(2,m_1,\ldots,m_r), \quad \text{$l$ even},
\end{equation}
\begin{equation}
\label{eq:epi_orient_holes_plus}
{\rm Epi}_o^+(\pi_1(O),\Z_{2l})=m^{2\mathfrak{g}+h-1}\cdot J_{2\mathfrak{g}+h-1}\Bigl(\dfrac{l}{m}\Bigr) \cdot \prod\limits_{i=1}^r\phi(m_i),\qquad m=\lcm(m_1,\ldots,m_r), \qquad \text{$l$ odd},
\end{equation}
where $J_k(n)$ is the Jordan's totient function. Here the number of epimorphisms is equal to zero if the argument $n$ of $J_k(n)$ is not an integer. 
\end{conseq_eng}

\evidpEng We prove the formula (\ref{eq:epi_orient_holes}), the formula (\ref{eq:epi_orient_holes_plus}) is proved similarly. From the Theorem \ref{theor:num_epimorph} and the Proposition \ref{propos:num_hom} we get
$$
\Epi_o(\pi_1(O), \Z_{l}) = \sum\limits_{d\mid l} \mu(l/d) \prod\limits_{i=1}^r \phi(m_i) d^{2{\mathfrak{g}}+h-1},
$$
where $2|d$ and all $m_i|d$. Define $m'$ as
$$
m':=\lcm(2,m_1,m_2,\ldots,m_r).
$$
Then
$$
\Epi_o(\pi_1(O), \Z_{l}) = \prod\limits_{i=1}^r \phi(m_i) \sum\limits_{m'\mid d\mid l} \mu(l/d) \cdot d^{2{\mathfrak{g}}+h-1} = 
\prod\limits_{i=1}^r \phi(m_i) \sum\limits_{d'\mid\frac{l}{m'}} \mu\Bigl(\dfrac{l}{m'd'}\Bigr)\cdot (m')^{2{\mathfrak{g}}+h-1} \cdot (d')^{2{\mathfrak{g}}+h-1}=
$$
$$
=(m')^{2{\mathfrak{g}}+h-1} \cdot J_{2{\mathfrak{g}}+h-1}\Bigl(\frac{l}{m'}\Bigr) \cdot \prod_{i=1}^r \phi(m_i),\qquad\qquad \text{$l$ even.}
$$
\qed

\begin{rem_eng}
It can be shown that for the case of a non-orientable orbifold $O$ with $h>0$ boundary components the formulas for the number of homomorphisms as well as formulas for the number of epimorphisms coincide with the formulas (\ref{eq:num_ord_pres_hom_orient}), (\ref{eq:num_ord_pres_hom_orient_plus}) and (\ref{eq:epi_orient_holes}), (\ref{eq:epi_orient_holes_plus}) with the difference that $2\mathfrak{g}$ is replaced by $\mathfrak{g}$. In particular, the numbers of epimorphisms are equal to
\begin{equation}
\label{eq:epi_non_orient_holes}
{\rm Epi}_o(\pi_1(O),\Z_{l})=(m')^{\mathfrak{g}+h-1}\cdot J_{\mathfrak{g}+h-1}\Bigl(\dfrac{l}{m'}\Bigr) \cdot \prod\limits_{i=1}^r\phi(m_i),\quad m'=\lcm(2,m_1,\ldots,m_r), \quad \text{$l$ even},
\end{equation}
\begin{equation}
\label{eq:epi_non_orient_holes_plus}
{\rm Epi}_o^+(\pi_1(O),\Z_{2l})=m^{\mathfrak{g}+h-1}\cdot J_{\mathfrak{g}+h-1}\Bigl(\dfrac{l}{m}\Bigr) \cdot \prod\limits_{i=1}^r\phi(m_i),\qquad m=\lcm(m_1,\ldots,m_r), \qquad \text{$l$ odd}.
\end{equation}
\end{rem_eng}

\begin{propos_eng}\label{propos:num_hom_nonorient_without_border}
The number of order-preserving homomorphisms from $\pi_1(O)$ to $\Z_d$ for the case of a non-orientable orbifold $O$ without boundary is equal to 
$$
\Hom_o(\pi_1(O),\Z_{d})=
\begin{cases}
d^{\mathfrak{g}-1} \prod\limits_{i=1}^r\phi(m_i), &\text{if $d$ is odd and $m_i|d$;}\\
2d^{\mathfrak{g}-1} \prod\limits_{i=1}^r\phi(m_i), &\text{if $d$ is even, $m_i|d$ and $\sum\limits_{i=1}^r d/m_i$ is odd;}\\
0 &\text{otherwise.}
\end{cases}
$$
\end{propos_eng}

\evidpEng For a non-orientable orbifold $O$ without boundary the relations for the fundamental group $\pi_1(O)$ 
$$
\prod\limits_{i=1}^rx_i\prod\limits_{k=1}^{\mathfrak{g}}a_k^2=1,\qquad\qquad x_i^{m_i}=1
$$
under the order-preserving homomorphism $\psi:\pi_1(O)\to \Z_d$ are transformed into equalities
\begin{equation}
\label{eq:psi_x_i_psi_a_k}
\sum\limits_{i=1}^r\psi(x_i)+2\sum\limits_{k=1}^{\mathfrak{g}} \psi(a_k)=0,\qquad\qquad 
m_i\cdot \psi(x_i)=0.
\end{equation}
Since the order of $\psi(x_i)$ is equal to $m_i$, from the second equality in (\ref{eq:psi_x_i_psi_a_k}) it follows that $m_i|d$ and there are $\phi(m_i)$ possible values of $\psi(x_i)$. 

Now from the first equality in (\ref{eq:psi_x_i_psi_a_k}) we can determine the number of possible values of $\psi(a_k)$. Rewrite this equality as
\begin{equation}
\label{eq:linear_congruence}
2z\equiv s\quad(\mod\,d),\qquad\qquad s:=d-\sum\limits_{i=1}^r\psi(x_i),\qquad \qquad z:=\sum\limits_{k=1}^{\mathfrak{g}} \psi(a_k).
\end{equation}
For odd $d$ the numbers $d$ and $2$ are coprime, so for a fixed $s$ there is only one possible value for $z$. It follows that there are $d^{\mathfrak{g}-1}$ possible values of $\psi(a_k)$, $k=1,\ldots,\mathfrak{g}$. 

For even $d$ the congruence relation (\ref{eq:linear_congruence}) can be satisfied only if $s$ is even. In this case there are exactly two possible values for $z$, so the number of homomorphisms is equal to
$$
2\cdot d^{\mathfrak{g}-1}\prod\limits_{i=1}^r\phi(m_i).
$$
It remains to determine the conditions for $s$ to be even. Since
$$
k=2k'=d-\sum\limits_{i=1}^r\psi(x_i)=2d'-\sum\limits_{i=1}^r\psi(x_i),
$$
the sum of all $\psi(x_i)$ must be even. We know that the order of $\psi(x_i)$ is equal to $m_i$. Hence
$$
\gcd(d,\psi(x_i))=\dfrac{d}{m_i}.
$$
Since $d$ is even, the parity of $\psi(x_i)$ coincides with the parity of $d/m_i$. As a consequence, in the case of even $d$ the number of homomorphisms is nonzero only if $\sum d/m_i$ is even. \qed

Analogous considerations can be used to prove the following proposition.

\begin{propos_eng}\label{propos:num_hom_nonorient_without_border_plus}
The number of orientation-and-order-preserving homomorphisms from $\pi_1(O)$ to $\Z_{2d}$ for the case of non-orientable orbifold $O$ without boundary is equal to 
$$
\Hom_o^+(\pi_1(O),\Z_{2d})=
\begin{cases}
d^{{\mathfrak{g}}-1}\prod\limits_{i=1}^r \phi(m_i),& \text{$d$ is odd and $m_i|d$;}\\
2d^{{\mathfrak{g}}-1}\prod\limits_{i=1}^r \phi(m_i),& \text{$d$ is even, $m_i|d$ and 
$\sum\limits_{i=1}^r \dfrac{d}{m_i}\equiv \mathfrak{g}\quad (\mod \,2)$;}\\
0&\text{otherwise.}
\end{cases}
$$
\end{propos_eng}

\begin{conseq_eng}
The number of order-preserving epi\-mor\-phisms for the case of a non-orientable orbifold $O$ without boundary in the case of odd $l$ is equal to
\begin{equation}
\label{eq:epi_nonorient_without_bound_odd}
\Epi_o(\pi_1(O), \Z_{l}) =m^{{\mathfrak{g}}-1} \cdot J_{{\mathfrak{g}}-1}\Bigl(\frac{l}{m}\Bigr) \cdot \prod_{i=1}^r \phi(m_i),\qquad m=\lcm(m_1,\ldots,m_r).
\end{equation}
For $l=2^qk$, where $k$ is odd and $q>1$, the number of epimorphisms is equal to 
\begin{equation}
\label{eq:epi_nonorient_without_bound_even_1}
\Epi_o(\pi_1(O), \Z_{l}) = 
2(m')^{\mathfrak{g}-1}\cdot J_{\mathfrak{g}-1}\Bigl(\dfrac{l}{m'}\Bigr) \cdot \prod\limits_{i=1}^r\phi(m_i),
\qquad m'=\lcm(2,b,m_1,\ldots,m_r),
\end{equation}
where $b$ is the denominator of the fraction $\sum\dfrac{1}{2m_i}$ after simplification. Finally, for $l=2k$ where $k$ is odd, we have
\begin{equation}
\label{eq:epi_nonorient_without_bound_even_2}
\Epi_o(\pi_1(O), \Z_{l}) = 
2(m')^{\mathfrak{g}-1}\cdot J_{\mathfrak{g}-1}\Bigl(\dfrac{l}{m'}\Bigr) \cdot \prod\limits_{i=1}^r\phi(m_i)
-m^{\mathfrak{g}-1}J_{\mathfrak{g}-1}\Bigl(\dfrac{l}{2m}\Bigr)\prod\limits_{i=1}^r\phi(m_i),
\end{equation}
where $m$ and $m'$ are defined as in the equations (\ref{eq:epi_nonorient_without_bound_odd}) and (\ref{eq:epi_nonorient_without_bound_even_1}).
\end{conseq_eng}

\begin{conseq_eng}
The numbers of orientation-and-order-preserving epi\-mor\-phisms for the case of a non-orientable orbifold $O$ without boundary are equal to
\begin{equation}
\label{eq:epi_nonorient_odd}
\Epi_o^+(\pi_1(O), \Z_{2l}) =m^{{\mathfrak{g}}-1} \cdot J_{{\mathfrak{g}}-1}\Bigl(\frac{l}{m}\Bigr) \cdot \prod_{i=1}^r \phi(m_i),\qquad m=\lcm(m_1,\ldots,m_r), \qquad\text{$l$ odd,}
\end{equation}
\begin{equation}
\label{eq:epi_nonorient_even}
\Epi_o^+(\pi_1(O), \Z_{2l}) = 2 \cdot (m')^{{\mathfrak{g}}-1} \cdot J_{{\mathfrak{g}}-1}\Bigl(\dfrac{l}{m'}\Bigr) \cdot \prod_{i=1}^r \phi(m_i),\quad m'=\lcm(2,m_1,\ldots,m_r), \quad  \text{$l$ even.}
\end{equation}
\end{conseq_eng}

{\setlength{\tabcolsep}{3pt}
\begin{table}[t!]
\begin{center}
\footnotesize
\begin{tabular}{ccccccc|cccccccc}
\midrule
$X$ & $l$ & $\pm$ & $\chi$ & $h$ & $m_i$ & $\Epi_o(\pi_1(O), \Z_{l})$ & $X$ & $l$ & $\pm$ & $\chi$ & $h$ & $m_i$ & $\Epi_o(\pi_1(O), \Z_{l})$ & $\Epi_o^+(...)$ \\ 
\midrule
$1$ & $1$ &$-$& $1$ & $0$ & [] & $1$ & $0$ & $1$ &$-$& $0$ & $0$ & [] & $1$ & $0$ \\
$1$ & $2$ & + & $1$ & $1$ & [$2$] & $1$ & $0$ & $2$ & + & $1$ & $1$ & [$2$ $2$] & $1$ & $0$ \\
$1$ & $3$ &$-$& $1$ & $0$ & [$3$] & $2$ & $0$ & $2$ &$-$& $1$ & $0$ & [$2$ $2$] & $2$ & $0$ \\
$1$ & $4$ & + & $1$ & $1$ & [$4$] & $2$ & $0$ & $2$ & + & $0$ & $2$ & [] & $2$ & $1$ \\
$1$ & $5$ &$-$& $1$ & $0$ & [$5$] & $4$ & $0$ & $2$ &$-$& $0$ & $0$ & [] & $3$ & $1$ \\
$1$ & $6$ & + & $1$ & $1$ & [$6$] & $2$ & $0$ & $2$ &$-$& $0$ & $1$ & [] & $2$ & $1$ \\
$1$ & $7$ &$-$& $1$ & $0$ & [$7$] & $6$ & $0$ & $3$ &$-$& $0$ & $0$ & [] & $2$ & $0$ \\
$1$ & $8$ & + & $1$ & $1$ & [$8$] & $4$ & $0$ & $4$ & + & $0$ & $2$ & [] & $2$ & $0$ \\
$1$ & $9$ &$-$& $1$ & $0$ & [$9$] & $6$ & $0$ & $4$ &$-$& $0$ & $0$ & [] & $4$ & $4$ \\
$1$ & $10$ & + & $1$ & $1$ & [$10$] & $4$ & $0$ & $4$ &$-$& $0$ & $1$ & [] & $2$ & $0$ \\
$1$ & $11$ &$-$& $1$ & $0$ & [$11$] & $10$ & $0$ & $5$ &$-$& $0$ & $0$ & [] & $4$ & $0$ \\
$1$ & $12$ & + & $1$ & $1$ & [$12$] & $4$ & $0$ & $6$ & + & $0$ & $2$ & [] & $4$ & $2$ \\ 
 ... & ... & ... & ... & ... & ... & ... & ... & ... & ... & ... & ... & ... & ... & ... \\
\midrule
$-1$ & $1$ &$-$& $-1$ & $0$ & [] & $1$ & $-2$ & $1$ &$-$& $-2$ & $0$ & [] & $1$ & $0$ \\
$-1$ & $2$ & + & $1$ & $1$ & [$2$ $2$ $2$] & $1$ & $-2$ & $2$ & + & $1$ & $1$ & [$2$ $2$ $2$ $2$] & $1$ & $0$ \\
$-1$ & $2$ & + & $0$ & $2$ & [$2$] & $2$ & $-2$ & $2$ &$-$& $1$ & $0$ & [$2$ $2$ $2$ $2$] & $2$ & $0$ \\
$-1$ & $2$ &$-$& $0$ & $1$ & [$2$] & $2$ & $-2$ & $2$ & + & $0$ & $2$ & [$2$ $2$] & $2$ & $0$ \\
$-1$ & $3$ &$-$& $1$ & $0$ & [$3$ $3$] & $4$ & $-2$ & $2$ &$-$& $0$ & $0$ & [$2$ $2$] & $4$ & $0$ \\
$-1$ & $4$ & + & $1$ & $1$ & [$2$ $4$] & $2$ & $-2$ & $2$ &$-$& $0$ & $1$ & [$2$ $2$] & $2$ & $0$ \\
$-1$ & $6$ & + & $1$ & $1$ & [$2$ $3$] & $2$ & $-2$ & $2$ & + & $-1$ & $1$ & [] & $4$ & $1$ \\
& &  & & & & & $-2$ & $2$ & + & $-1$ & $3$ & [] & $4$ & $1$ \\
& &  & & & & & $-2$ & $2$ &$-$& $-1$ & $0$ & [] & $7$ & $1$ \\
& &  & & & & & $-2$ & $2$ &$-$& $-1$ & $1$ & [] & $4$ & $1$ \\
& &  & & & & & $-2$ & $2$ &$-$& $-1$ & $2$ & [] & $4$ & $1$ \\
& &  & & & & & $-2$ & $3$ &$-$& $0$ & $0$ & [$3$] & $6$ & $0$ \\
& &  & & & & & $-2$ & $4$ &$-$& $1$ & $0$ & [$2$ $2$ $2$] & $2$ & $2$ \\
& &  & & & & & $-2$ & $4$ & + & $1$ & $1$ & [$4$ $4$] & $4$ & $0$ \\
& &  & & & & & $-2$ & $4$ &$-$& $1$ & $0$ & [$4$ $4$] & $8$ & $0$ \\
& &  & & & & & $-2$ & $4$ & + & $0$ & $2$ & [$2$] & $2$ & $0$ \\
& &  & & & & & $-2$ & $4$ &$-$& $0$ & $0$ & [$2$] & $8$ & $0$ \\
& &  & & & & & $-2$ & $4$ &$-$& $0$ & $1$ & [$2$] & $2$ & $0$ \\
& &  & & & & & $-2$ & $6$ & + & $1$ & $1$ & [$2$ $6$] & $2$ & $0$ \\
& &  & & & & & $-2$ & $6$ &$-$& $1$ & $0$ & [$2$ $6$] & $4$ & $0$ \\
& &  & & & & & $-2$ & $6$ & + & $1$ & $1$ & [$3$ $3$] & $4$ & $4$ \\
& &  & & & & & $-2$ & $6$ &$-$& $1$ & $0$ & [$3$ $3$] & $4$ & $4$ \\
& &  & & & & & $-2$ & $8$ &$-$& $1$ & $0$ & [$2$ $4$] & $4$ & $4$ \\
& &  & & & & & $-2$ & $12$ &$-$& $1$ & $0$ & [$2$ $3$] & $4$ & $4$ \\
\midrule
\end{tabular}
\caption{Orbifolds covered by Euler characteristic $X$ surfaces. For $X \geq 0$ only the first few rows are listed. For even $X$ all $\Epi_o(\pi_1(O), \Z_{l})$ epimorphisms are split into those corresponding to coverings by all and only by orientable surfaces of Euler characteristic $X$.}
\label{table:all_orbifolds}
\end{center}
\end{table}}

The results of calculations using the formulas (\ref{eq:epi_orient_holes})--(\ref{eq:epi_non_orient_holes_plus}) and (\ref{eq:epi_nonorient_without_bound_odd})--(\ref{eq:epi_nonorient_even}) are given in the Table \ref{table:all_orbifolds}. In the first column of this table we give the Euler characteristic $X$ of the covering surface, in the second column --- the period $l$ of the homeomorphism $h$. A plus sign in the third column means that the covered surface is orientable, and a minus sign — that it is non-orientable. The fourth column shows the Euler characteristic $\chi$ of the covered surface. In the fifth column we give the number of boundary components, and in the sixth column we give the list of branch point indices.

In the case of odd values of $X$ (the left side of the table), the covering surface is always non-orientable, so any epimorphism does not preserve the sign. The number of such epimorphisms is given in the seventh column of the table. In the case of even $X$ (the right side of the table), the covering surface can be either non-orientable or orientable. Order-preserving epimorphisms (the seventh column of the table) correspond to all surfaces of the Euler characteristic $X$. The difference between these numbers corresponds to coverings of orbifolds by non-orientable surfaces. 

In addition, the authors have implemented two programs which can extend this table up to a given genus $g$. These programs can be found at \url{https://github.com/krasko/o_r_orbifolds} and \url{https://github.com/krasko/n_o_orbifolds}.

\section{Counting quotient maps on the disk and on the Möbius band} 

The theorems \ref{theor:unsensed_orient} and \ref{theor:unsensed_nonorient} reduce the problem of enumerating unsensed maps to the problem of determining the numbers $\tau_O(n)$ of quotient maps on an orbifold $O$. As we mentioned before, in the case of unsensed maps such orbifolds may contain boundary components. In its turn, these boundary components may contain vertices and/or edges. This property significantly complicates enumeration of quotient maps on orbifolds. In this section we demonstrate the technique of counting such maps with the example of maps on the simplest orientable surface with a boundary, namely a disk, and with the example of maps on the simplest non-orientable surface with a boundary, namely the Möbius band.

\begin{figure}[ht]
\centering
	\centering
    	\includegraphics[scale=0.7]{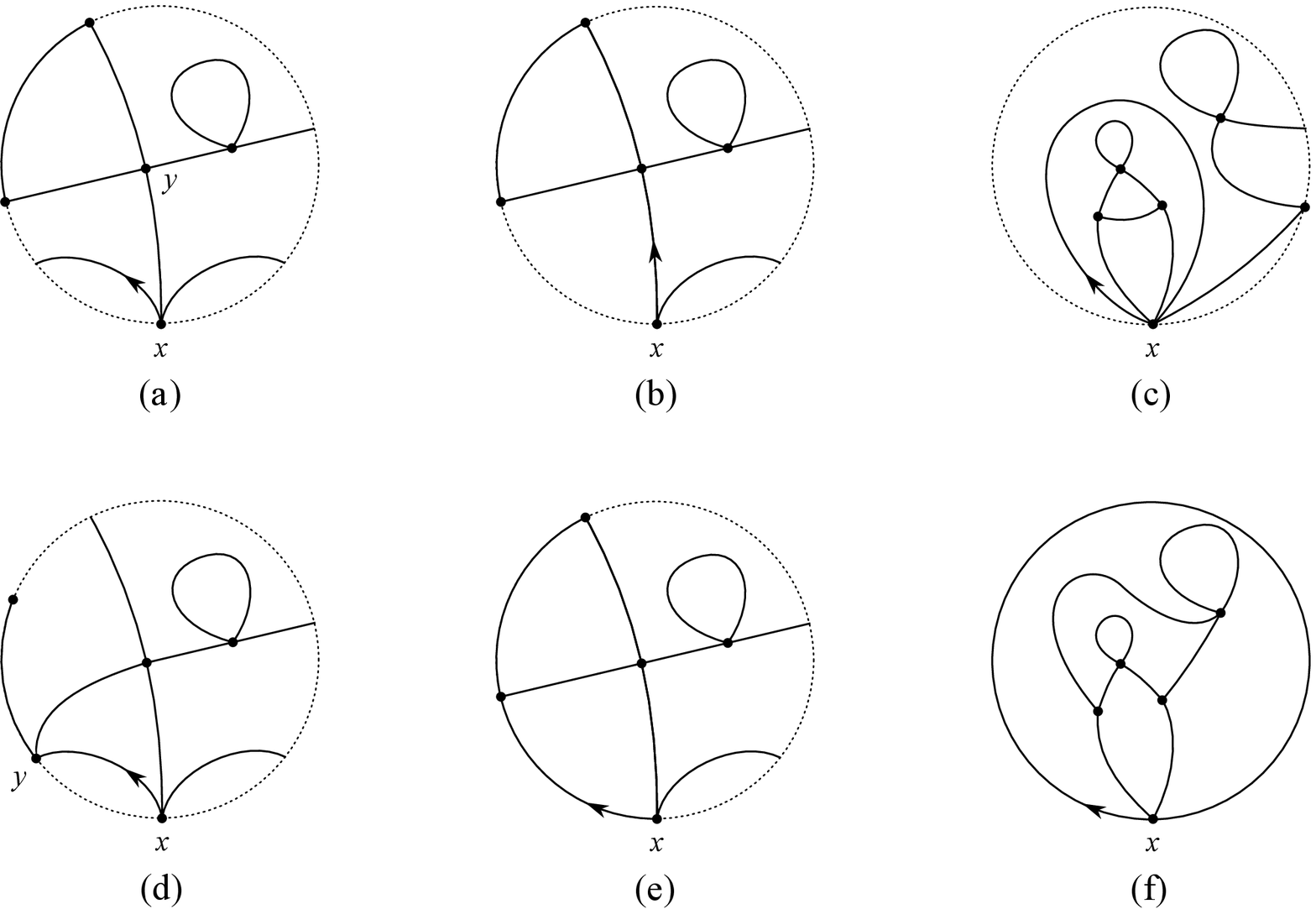}
	\caption{}
\label{fig:disk_all_maps}
\end{figure}

The recurrence relations for the numbers of maps will be different depending on the location of the root dart in these maps. We begin with the case when the root dart is incident to a vertex on the boundary, but does not lie on the boundary itself (see Figure \ref{fig:disk_all_maps}). Without loss of generality, we will assume that the root dart is the leftmost dart incident to the root vertex.

\begin{propos_eng}
The number $d^{(0)}_{n,k}$ of quotient maps on a disk in the case when the root dart doesn't lie on the boundary but is incident to a vertex $x$ of degree $k$ which lies on the boundary, can be expressed by the following formula:
$$
d^{(0)}_{n,k}=d^{(0)}_{n-1,k-1}+\sum\limits_{k'=0}^{n-k-1}d^{(0)}_{n-2,k-1+k'}+
\sum\limits_{n'=0}^{n-2}\sum\limits_{k'=0}^{k-2}d^{(0)}_{n',k'}\cdot s_{n-n'-2,k-k'-2}+
$$
\begin{equation}
\label{eq:d_n_k}
+\sum\limits_{n'=0}^{n-2}\sum\limits_{k'=0}^{n-k-1}\sum\limits_{k''=0}^{k'}\sum\limits_{i_1=0}^{1}\sum\limits_{i_2=0}^{1}
d^{(i_1)}_{n-2-n',k+k''-1}\cdot d^{(i_2)}_{n',k'-k''}.
\end{equation}
Here $n$ is the number of darts; $d^{(1)}_{n,k}$ is the number of quotient maps on a disc with the root dart that belongs to a boundary edge, $k$ other darts incident to the root which all lie in the interior of the disc and $n$ darts in total; $s_{n, k}$ is the number of maps with $n$ darts and the root vertex of degree $k$ on a sphere.
\end{propos_eng}

\evidpEng Indeed, the first summand in (\ref{eq:d_n_k}) corresponds to the case of the root dart being a halfedge (Figure \ref{fig:disk_all_maps}(a)). Contracting this dart yields a map with the number of darts and the degree of the root vertex decreased by 1. The second summand in (\ref{eq:d_n_k}) corresponds to contracting a root edge that joins the root vertex $x$ and an internal vertex $y$ (Figure \ref{fig:disk_all_maps}(b)). The number of darts is decreased by two in this case, and the new degree of the root vertex becomes equal to the sum of those of $x$ and $y$ minus two. In the third case the root edge is a loop (Figure \ref{fig:disk_all_maps}(c)). Contracting this loop splits the map into two maps, one on a disc and one on a sphere. The double summation iterates over all possible degrees of the new root vertices and the distribution of darts among the two maps. Finally, the last term in (\ref{eq:d_n_k}) corresponds to contracting an edge that connects the root vertex $x$ with another vertex $y$ on the boundary (Figure \ref{fig:disk_all_maps}(d)). In the general case, the non-root vertex $y$ may have $(1)$ or not have $(0)$ a boundary edge on either side. Consequently, we have four different options described by the indices $i_1$ and $i_2$. The index $k''$ enumerates darts incident to $y$ that lie to the left of the root edge. The index $k'$ enumerates all darts incident to $y$ except for possible boundary edges and the dart that belongs to the root edge.  \qed

Now we consider the case of maps in which the root dart is the boundary dart incident to the vertex $x$ lying on the boundary (see Figure \ref{fig:disk_all_maps}(e)-(f)).

\begin{propos_eng}
The number $d^{(1)}_{n,k}$ of maps in which the root dart is the only boundary dart incident to $x$ can be counted as
$$
d^{(1)}_{n,k}=\sum\limits_{k'=0}^{n-k-2}\left(d^{(0)}_{n-1,k+k'}+d^{(1)}_{n-1,k+k'}\right),
$$ 
and the number $d^{(2)}_{n,k}$ of rooted quotient maps that have two boundary darts incident to the root vertex $x$ is calculated by the formula
$$
d^{(2)}_{n,k}=s_{n-1,k}+\sum\limits_{k'=0}^{n-k-2}\left(d^{(0)}_{n-1,k+k'}+d^{(1)}_{n-1,k+k'}\right).
$$
\end{propos_eng}

\evidpEng Indeed, contracting the root boundary edge $\{x,y\}$ yields, depending on the number of boundary edges incident to $y$, a map with $n-1$ darts with ($d^{(1)}_{n-1,k+k'}$) or without ($d^{(0)}_{n-1,k+k'}$) a boundary edge incident to the root. The index $k'$ stands for the number of darts incident to $y$ lying in the interior of the disc. 

The summand $s_{n-1,k}$ in the formula for $d^{(2)}_{n,k}$ corresponds to the case when the root boundary edge is a loop (see Figure \ref{fig:disk_all_maps}(f)): contracting this loops yields a map on a sphere with  $n-1$ darts and with the root vertex of degree $k$. \qed

It remains to consider the case of rooted quotient maps with the root vertex lying in the interior of the disk.

\begin{propos_eng} 
The number $d_{n,k}$ of rooted quotient maps with the root vertex $x$ lying in the interior of the disk can be counted using the following recurrence relation: 
$$
d_{n,k}=d^{(0)}_{n-1,k-1}+\sum\limits_{k'=0}^{n-k-1}d_{n-2,k+k'-1}+
2\sum\limits_{n'=0}^{n-2}\sum\limits_{k'=0}^{k-2}d_{n',k'}\cdot s_{n-n'-2,k-k'-2}+
$$
$$
+\sum\limits_{k'=0}^{n-k-1}\sum\limits_{i_1=0}^{1}\sum\limits_{i_2=0}^{1}
(k'+1)\cdot d^{(i_1+i_2)}_{n-2,k+k'-1}.
$$
\end{propos_eng}

\evidpEng The proof of this formula almost completely repeats the derivation of formula (\ref{eq:d_n_k}). The multiplier $2$ in the third summand appears here because the root edge that is a loop can be oriented in two possible ways. The final summand corresponds to contracting an edge $\{x,y\}$ where $y$ lies on the boundary. The multiplier $k'+1$ denotes the number of non-boundary darts incident to $y$. Since we assume that the root dart incident to the new root vertex will always be the leftmost dart, $k'+1$ different quotient maps will be reduced to the same map with the root vertex lying on the boundary. \qed

\begin{figure}[ht]
\centering
	\centering
    	\includegraphics[scale=0.7]{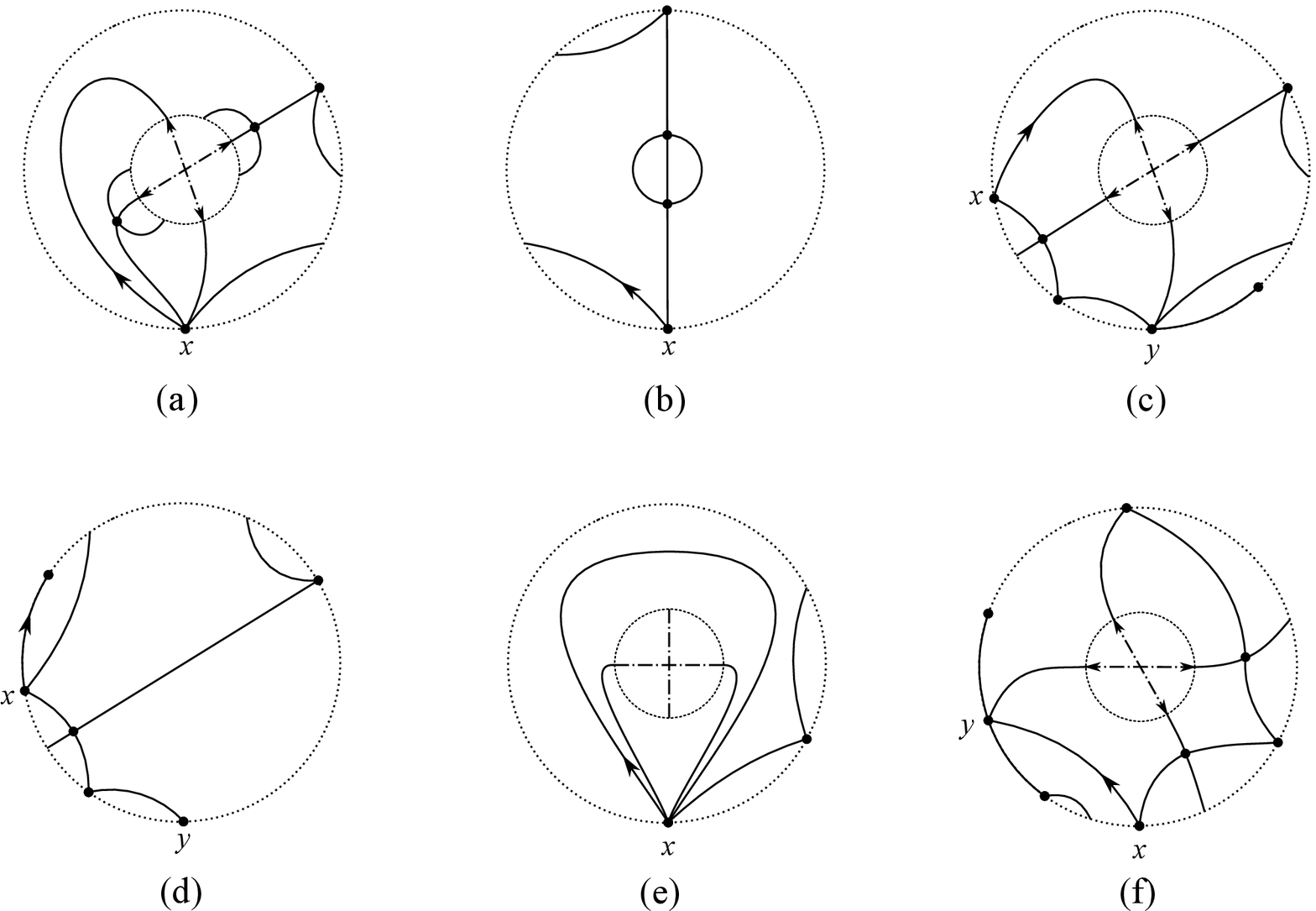}
	\caption{}
\label{fig:moebius}
\end{figure}

Now let's turn to the M\"obius band. We will consider in detail the case of maps having the root vertex on the boundary and having no boundary edges incident to the root vertex (Figure \ref{fig:moebius}; it will be convenient to use the representation of a M\"obius band as an annulus with a cross-cap glued to its inner boundary). 

\begin{propos_eng} 
The recurrence relation for the numbers $m^{(0)}_{n,k}$ of maps with the root vertex $x$ lying on the boundary and having no boundary edges incident to the root vertex has the following form:
$$
m^{(0)}_{n,k} = m^{(0)}_{n-1,k-1} + \sum\limits_{k'=0}^{n-k-1}m^{(0)}_{n-2,k-1+k'} + (k-1) \cdot d^{(0)}_{n-2,k-2}+ \sum\limits_{k'=0}^{n-k-1}\sum\limits_{i_1=0}^{1}\sum\limits_{i_2=0}^{1}d^{(i_1, \{0, i_2\})}_{n-2,k-1+k'}+
$$
\begin{equation}
\label{eq:m_n_k}
+\sum\limits_{n'=0}^{n-2}\sum\limits_{k'=0}^{k-2}(d^{(0)}_{n',k'}\cdot p_{n-n'-2,k-k'-2}+m^{(0)}_{n',k'}\cdot s_{n-n'-2,k-k'-2})+
\end{equation}
$$
+ \sum\limits_{n'=0}^{n-2}\sum\limits_{k'=0}^{n-k-1}\sum\limits_{k''=0}^{k'}\sum\limits_{i_1=0}^{1}\sum\limits_{i_2=0}^{1}
(m^{(i_1)}_{n-2-n',k+k''-1}\cdot d^{(i_2)}_{n',k'-k''} + d^{(i_1)}_{n-2-n',k+k''-1}\cdot m^{(i_2)}_{n',k'-k''}).
$$
\end{propos_eng}

\evidpEng Indeed, the first two summands are analogous to the corresponding summands in (\ref{eq:d_n_k}) and describe the cases of the root dart being a half-edge and the root edge joining the root vertex with a vertex in the interior of the surface. The third summand in the right-hand side of (\ref{eq:m_n_k}) describes maps which have the root edge going through the cross-cap (Figure \ref{fig:moebius}(a)). In this case cutting along this edge, inverting one of the obtained parts and gluing these parts back yields a quotient map on a disc (Figure \ref{fig:moebius}(b)). The multiplier $(d-1)$ has the same combinatorial meaning as the analogous multiplier for maps on the projective plane (see \cite{Bender_1988}).

The fourth summand in the right-hand side of (\ref{eq:m_n_k}) describes the situation when the root edge goes thought the cross-cap and joins the root vertex $x$, incident to $k$ internal darts, with a different vertex $y$ lying on the boundary (Figure  \ref{fig:moebius}(c)). The vertex $y$ may have or not have a boundary edge both on its right ($i_1=1$ or $i_1=0$) and on its left side ($i_2=1$ or $i_2=0$). It can also have some amount $k'$ of interior darts lying to the right of the root edge $\{x,y\}$. So, for the map depicted in Figure \ref{fig:moebius}(c), $i_1=1$, $i_2=0$, and $k'=1$. Contracting the root edge yields a map on a disc with a new root vertex $x'$ and an additional distinguished vertex $y'$ on its boundary (Figure \ref{fig:moebius}(d)). All darts incident to $y$ and lying on the right of the root edge will be incident to the new root vertex $x'$, whereas darts that lied to the left of it will be incident to the distinguished vertex $y'$. Note that after this operation the vertex $y'$ cannot have an incident boundary edge on its right side. The number of maps with such properties in the formula (\ref{eq:m_n_k}) is denoted by $d^{(i_1, \{0, i_2\})}_{n-2,k-1+k'}$.

The next summand in (\ref{eq:m_n_k}) enumerates maps with the root edge being a loop that doesn't go through the cross-cap. If the crosscap lies inside this loop, contracting the loop yields a map on a disc and a map on a projective plane (Figure \ref{fig:moebius}(e)). If the loop does not enclose the crosscap, contracting it yields a sphere and a M\"obius band. These cases correspond to the summands $d^{(0)}_{n',k'}\cdot p_{n-n'-2,k-k'-2}$ and $m^{(0)}_{n',k'}\cdot s_{n-n'-2,k-k'-2}$ in (\ref{eq:m_n_k}), where $p_{n, k}$ is the number of projective plane maps with $k$ darts incident to the root and $n$ darts in total. 

Finally, the last two summands in the right-hand side of (\ref{eq:m_n_k}) describe the situation when the root edge joins the root vertex $x$ with another vertex $y$ on the boundary and does not pass through the cross-cap (Figure \ref{fig:moebius}(f)). Depending on the side of the root edge which the cross-cap lies in, we obtain one or another summand in (\ref{eq:m_n_k}). Here the numbers $m^{(1)}_{n,k}$ enumerate maps on a M\"obius band which are defined analogously to the maps on a disc enumerated by $d_{n,k}^{(1)}$. \qed

\section{Counting quotient maps on an arbitrary orbifold}

A disc is the simplest example of an orientable surface with a boundary. In order to use the approach of Tutte in the general case of a sphere with $g$ handles and $h$ boundary components, we will need to overcome a number of technical difficulties. To conclude what parameters the recurrence relations should depend on, next we describe some of these difficulties.

\begin{figure}[ht]
\centering
	\centering
    	\includegraphics[scale=0.6]{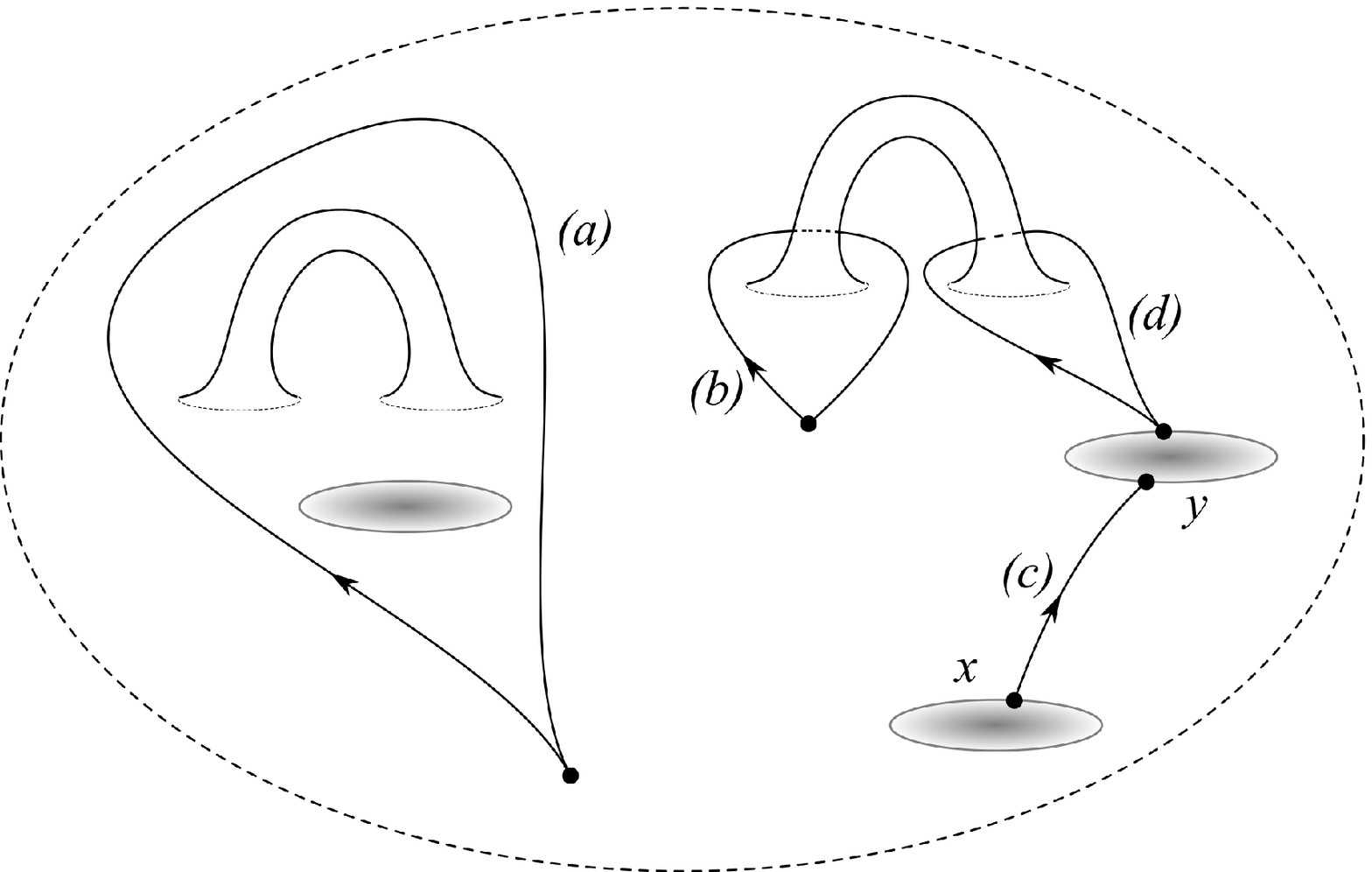}
	\caption{}
\label{fig:topology_1}
\end{figure}

First assume that a map contains a loop that wraps $h_1$ boundary components and $g_1$ handles (fig. \ref{fig:topology_1}(a)). Contracting this loop yields two surfaces with maps on them. We've already encountered a special case of this situation for the disc: contracting a loop yields a new disc and a sphere. Consequently, besides the numbers $n$ of darts and the degree $k$ of the root vertex, we also need to keep track of the numbers $h$ of boundary components and $g$ of handles. 

Now assume that the root edge is a loop wrapping a handle (Figure \ref{fig:topology_1}(b)). Contracting this loop splits the root vertex of degree $k$ into two vertices with the sum of degrees equal to $k-2$. For example, contracting a torus along its meridian yields a sphere with two distinguished vertices. In the general case, multiple applications of this procedure may yield $m$ internal distinguished vertices and we have to track their degrees. 

The third difficulty stems from the possibility of joining two boundary components or splitting a boundary component into two. The former case takes place when the root vertex connects two vertices on two different boundary components (Figure \ref{fig:topology_1}(c)). Contracting it merges these boundary components. In addition to the root vertex, a new distinguished vertex appears on the boundary, and we have to keep track of boundary edges incident to it. The latter case happens when the root edge wraps a handle, both vertices incident to it lie on the same boundary and they are distinct (Figure \ref{fig:topology_1}(d)). After contracting the root edge there appears a new distinguished vertex which lies on a separate boundary component. Consequently, we have to keep track of the list of distinguished vertices on each boundary component; each such vertex will be characterized by the presence or absence of boundary darts on each side.

The final difficulty is the possible presence of branch points on orbifolds. For surfaces without a boundary, branch points may fall into vertices, faces and dangling semiedges of quotient maps. It is possible to get rid of branch points in dangling semiedges using a combinatorial approach of double counting \cite{Mednykh_Nedela}. Using the Euler-Poincare formula one can then determine the total number of faces and vertices in a quotient map and reduce the problem of enumerating quotient maps to the problem of enumerating maps on surfaces without branch points.

As shown in \cite{Azevedo}, orbifolds with a boundary arising as quotients of orientable surfaces can't have dangling semiedges ending in branch points. However, the presence of a boundary makes it impossible to uniquely determine the number of faces and vertices which can coincide with branch points. Consequently, in this case we have to introduce one more parameter that would track the overall number of vertices and faces in the interior of the surface. 

So, the recurrence relations for the numbers of quotient maps on orientable surfaces with $h$ boundary components and $g$ handles will also depend on the number $n$ of darts, the degree $k$ of the root vertex, the list of $m$ degrees of internal vertices, and also on the lists of distinguished vertices for each of $h$ boundary components together with the information about boundary darts for each such vertex. If $h > 0$, they will also depend on the number of interior faces and vertices.

\begin{figure}[ht]
\centering
	\centering
    	\includegraphics[scale=0.6]{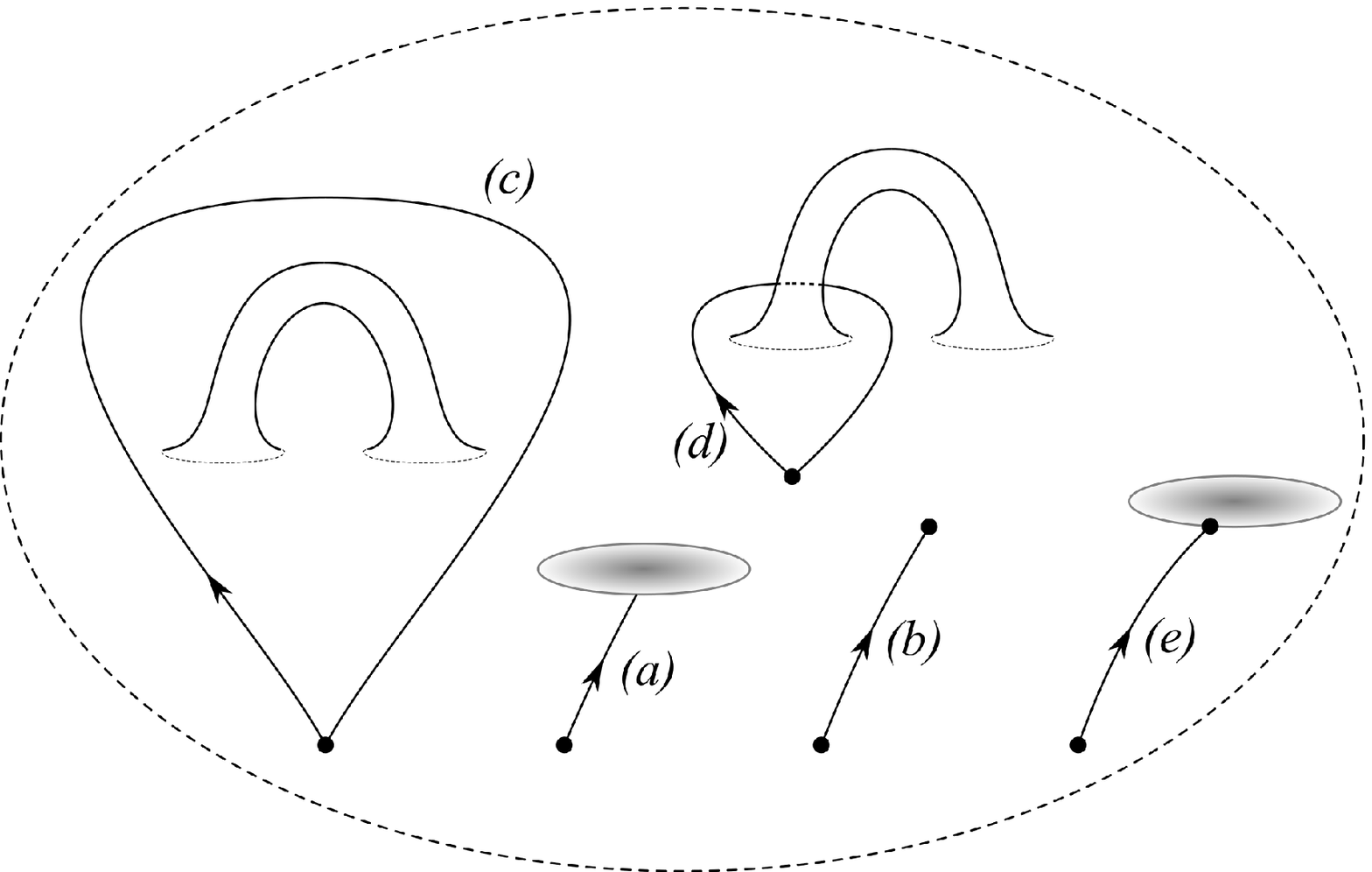}
	\caption{}
\label{fig:topology_2}
\end{figure}

Next we will describe a schema for obtaining these recurrence relations. Conceptually, it will be sufficient to consider two recurrence relations differing in the position of the root vertex. The first relation corresponds to the case of the root vertex $x$ located in the interior of a surface. We have to consider five different situations:
\begin{itemize}
\item[(a)] The root dart is a half edge, that is, it ends on a boundary (Figure \ref{fig:topology_2}(a));\\
\item[(b)] The root edge ends in an interior vertex $y$ different from $x$ (Figure \ref{fig:topology_2}(b));\\
\item[(c)] The root edge is a loop; contracting it splits the surface into two (Figure \ref{fig:topology_2}(c));\\
\item[(d)] The root edge is a loop; it wraps a handle so that contracting it decreases the surface genus (Figure \ref{fig:topology_2}(d));\\
\item[(e)] The root edge connects the root vertex $x$ with a vertex $y$ on a boundary (Figure \ref{fig:topology_2}(e)).
\end{itemize}

\begin{figure}[ht]
\centering
	\centering
    	\includegraphics[scale=0.6]{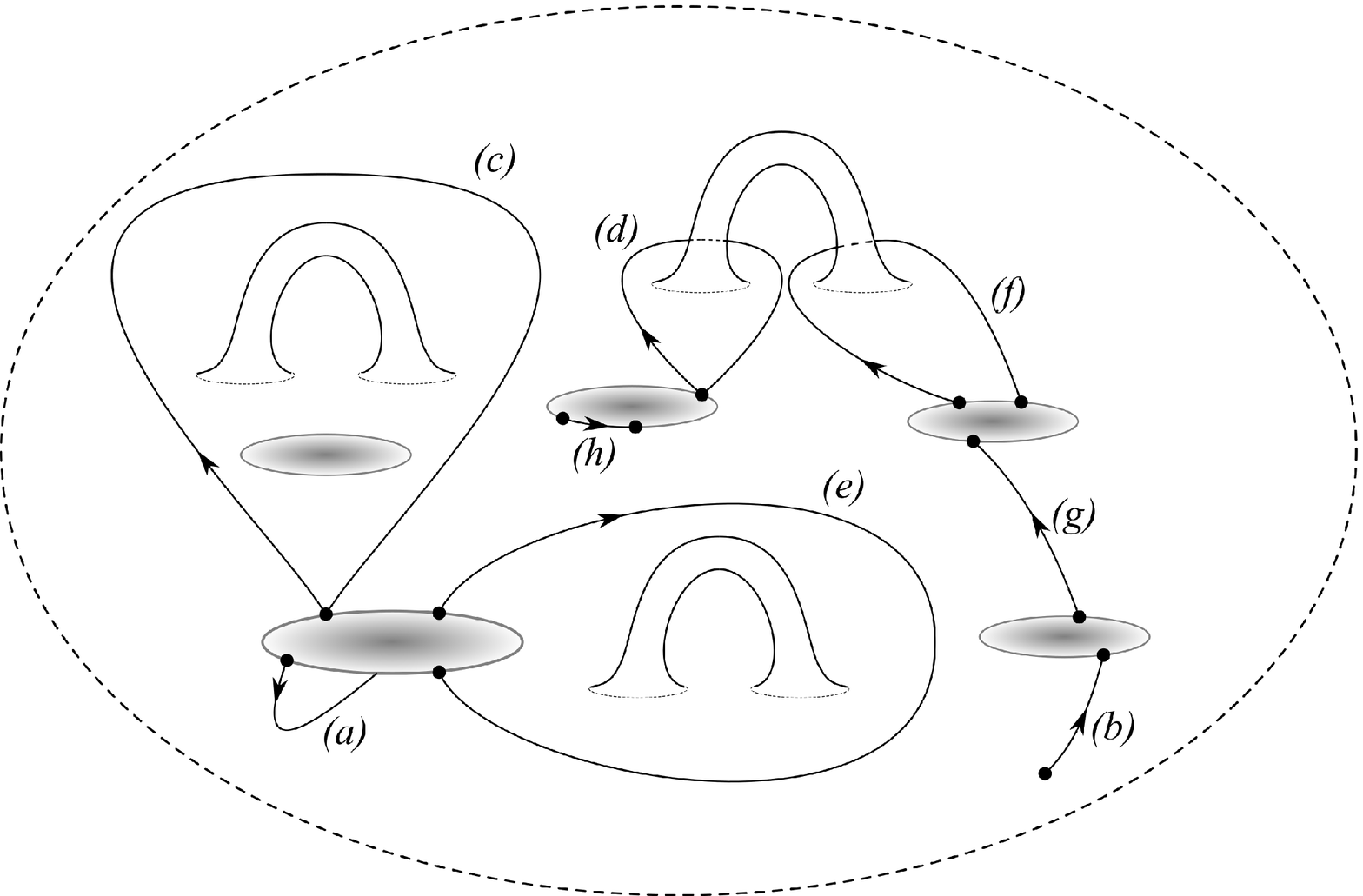}
	\caption{}
\label{fig:topology_3}
\end{figure}

The second recurrence relation enumerates maps with the root vertex $x$ lying on a boundary (Figure \ref{fig:topology_3}). The following topologically distinct possibilities may take place:
\begin{itemize}
\item[(a)] The root dart is a half edge ending on the boundary component that contains $x$ (Figure \ref{fig:topology_3}(a));\\
\item[(b)] The root edge ends in an interior vertex $y$ different from $x$ (Figure \ref{fig:topology_3}(b));\\
\item[(c)] The root edge is a loop; contracting it splits the surface into two (Figure \ref{fig:topology_3}(c));\\
\item[(d)] The root edge is a loop; it wraps a handle so that contracting it decreases the surface genus (Figure \ref{fig:topology_3}(d));\\
\item[(e)] The root edge connects the vertex $x$ with a different vertex $y$ on the same boundary component in a way that contracting the root edge $\{x, y\}$ splits the surface into two (Figure \ref{fig:topology_3}(e));\\
\item[(f)] The root edge connects the vertex $x$ with a different vertex $y$ on the same boundary component in a way that contracting the root edge $\{x, y\}$ decreases the surface genus (Figure \ref{fig:topology_3}(f));\\
\item[(g)] The root edge connects vertices on two different boundary components (Figure \ref{fig:topology_3}(g));\\
\item[(h)] The root edge is a boundary edge (Figure \ref{fig:topology_3}(h)); it could be a loop or not.
\end{itemize}

Now let's turn to non-orientable orbifolds. Any such orbifold is a sphere with $l + h$ boundary components, $l$ of which ($l > 0$) are glued with cross-caps. As the example of a M\"obius band shows, for non-orientable surface there is also an additional possibility that the root edge would go through a crosscap (Figure \ref{fig:moebius}(a) and \ref{fig:moebius}(c)). Consequently, instead of five cases of contracting the root edge for the root vertex $x$ in the interior and eight cases for the root vertex $x$ on the boundary, for non-orientable surfaces we will have six and ten cases, correspondingly.

It's also important to note that for the case of maps on non-orientable surfaces it is convenient to derive recurrence relations which enumerate maps regardless of orientability, that is, enumerate maps on orientable and non-orientable surfaces together \cite{Gao_1991}. The number of maps on non-orientable surfaces can then be obtained by subtracting the number of maps on orientable surfaces from this number.

Rewriting all these results in a form of a single mathematical formula is impractical, as it would be extremely cumbersome. Instead of that, authors have implemented two programs (for orientable and for non-orientable surfaces) that comprise all technical details of the above-described considerations. They can be used to calculate the numbers of quotient maps on orbifolds of a given signature and, provided a list of orbifolds, enumerate unsensed maps on surfaces of a given genus up to a certain number of edges. These programs can be found at \url{https://github.com/krasko/reflexible_maps} and \url{https://github.com/krasko/non-orientable_maps}.

\section*{Conclusion}

The results presented in this article allowed us to enumerate unsensed maps with $n$ edges on orientable and non-orientable surfaces of a given genus $g$. In the tables \ref{table:table_last_1_or}--\ref{table:table_last_3_or} we give the results of these computations for unsensed maps on orientable surfaces of genus $g\in[1,10]$ with $n\in[2,20]$ edges. Note that in  \cite{Walsh_generation} these numbers were obtained by explicit generation of the corresponding maps for $g\in[1,5]$ and $n\in[2,11]$. The method described in this work allows to advance further in unsensed maps enumeration. In the tables \ref{table:table_last_1}--\ref{table:table_last_3} we give the results of computations of unsensed maps on non-orientable surfaces of genus $g\in[1,13]$, for $n\in[1,14]$.

\begin{table}[h!]
\begin{center}
\footnotesize
\begin{tabular}{c|cccc}
\midrule
$ n\backslash g $ &\phantom{00000}$1$\phantom{00000}&\phantom{0000}$2$\phantom{0000}&\phantom{0000000000}$3$\phantom{0000000000}&\phantom{0000000}$4$\phantom{0000000}\\ 
\midrule
 2 &                       1 &                       0 &                       0 &                       0\\
 3 &                       6 &                       0 &                       0 &                       0\\
 4 &                      40 &                       4 &                       0 &                       0\\
 5 &                     320 &                      76 &                       0 &                       0\\
 6 &                    2946 &                    1395 &                      82 &                       0\\
 7 &                   29364 &                   24950 &                    4348 &                       0\\
 8 &                  309558 &                  427336 &                  160050 &                    7258\\
 9 &                 3365108 &                 6987100 &                 4696504 &                  688976\\
10 &                37246245 &               109761827 &               118353618 &                37466297\\
11 &               416751008 &              1668376886 &              2675297588 &              1512650776\\
12 &              4696232371 &             24689351504 &             55758114082 &             50355225387\\
13 &             53186743416 &            357467967214 &           1091344752470 &           1461269893538\\
14 &            604690121555 &           5083309341304 &          20318440463052 &          38236656513725\\
15 &           6896534910612 &          71209097157108 &         363171011546210 &         922552326544030\\
16 &          78867385697513 &         984963603696282 &        6275111078422480 &       20847359639841664\\
17 &         904046279771682 &       13477371260785608 &      105369657960443204 &      446290728182323620\\
18 &       10384916465797240 &      182698708325667710 &     1726590417107274316 &     9129236228868478458\\
19 &      119522063788612992 &     2456600457435363198 &    27699670730854989616 &   179639607187998993180\\
20 &     1378014272286250059 &    32796863046711248526 &   436246336648672487876 &  3418366706444416598777\\\midrule
\end{tabular}
\caption{Orientable genus $g$ maps with $n$ edges}
\label{table:table_last_1_or}
\end{center}
\end{table}

\begin{table}[h!]
\begin{center}
\footnotesize
\begin{tabular}{c|ccc}
\midrule
$ n\backslash g $ &\phantom{00000}$5$\phantom{00000}&\phantom{0000}$6$\phantom{0000}&\phantom{0000000000}$7$\phantom{0000000000}\\ 
\midrule
10 &                 1491629 &                       0 &                       0\\
11 &               195728778 &                       0 &                       0\\
12 &             14019733828 &               506855279 &                       0\\
13 &            724646387874 &             84930743344 &                       0\\
14 &          30220873171570 &           7601322881752 &            254118439668\\
15 &        1079253898643492 &         482475325333252 &          52148049818864\\
16 &       34231899372185491 &       24347701836204379 &        5634797561708385\\
17 &      988157793188200998 &     1038820801135250668 &      426497331688178676\\
18 &    26412878913430197293 &    38928478953655850016 &    25388940147173859412\\
19 &   662133032168309300424 &  1314623638623845390906 &  1265623233919838264624\\
20 & 15719783014093104131694 & 40749347642026348171659 & 54940200059090328012148\\
\midrule
\end{tabular}
\caption{Orientable genus $g$ maps with $n$ edges}
\label{table:table_last_2_or}
\end{center}
\end{table}

\begin{table}[h!]
\begin{center}
\footnotesize
\begin{tabular}{c|ccc}
\midrule
$ n\backslash g $ &\phantom{00000}$8$\phantom{00000}&\phantom{0000}$9$\phantom{0000}&\phantom{0000000000}$10$\phantom{0000000000}\\ 
\midrule
16 &         176377605783906 &                       0 &                       0\\
17 &       43058445711817178 &                       0 &                       0\\
18 &     5477393987229533288 &      162019808170348933 &                       0\\
19 &   483573171728920541590 &    46037869248765236030 &                       0\\
20 & 33299663456795126129156 &  6762460437287955976080 &   190375587419231088550\\
\end{tabular}
\caption{Orientable genus $g$ maps with $n$ edges}
\label{table:table_last_3_or}
\end{center}
\end{table}

\begin{table}[h!]
\begin{center}
\footnotesize
\begin{tabular}{c|ccccc}
\midrule
$ n\backslash g $ &\phantom{0000000}$1$\phantom{0000000}&\phantom{0000000}$2$\phantom{0000000}&\phantom{0000000}$3$\phantom{0000000}&\phantom{0000000}$4$\phantom{0000000}&\phantom{0000000}$5$\phantom{0000000}\\ 
\midrule
 1 &                       1 &                       0 &                       0 &                       0 &                       0\\
 2 &                       4 &                       2 &                       0 &                       0 &                       0\\
 3 &                      19 &                      16 &                       8 &                       0 &                       0\\
 4 &                     106 &                     137 &                     128 &                      47 &                       0\\
 5 &                     709 &                    1254 &                    1890 &                    1372 &                     473\\
 6 &                    5356 &                   12597 &                   27036 &                   31007 &                   22556\\
 7 &                   44558 &                  133518 &                  379491 &                  611322 &                  704066\\
 8 &                  397146 &                 1464725 &                 5229092 &                11017122 &                17691240\\
 9 &                 3716039 &                16373700 &                70805740 &               186044902 &               387965547\\
10 &                35967272 &               185086459 &               944106760 &              2992773591 &              7743850792\\
11 &               356784008 &              2106417804 &             12426068215 &             46378655568 &            144265686318\\
12 &              3605014966 &             24081813881 &            161793730426 &            697928466684 &           2549632137634\\
13 &             36955965852 &            276231542440 &           2087922762430 &          10257901689164 &          43225637985830\\
14 &            383320824698 &           3176840637522 &          26745380615078 &         147883893647230 &         708619408979790\\\midrule
\end{tabular}
\caption{Non-orientable genus $g$ maps with $n$ edges}
\label{table:table_last_1}
\end{center}
\end{table}

\begin{table}[h!]
\begin{center}
\footnotesize
\begin{tabular}{c|cccc}
\midrule
$ n\backslash g $ &\phantom{00000000}$6$\phantom{00000000}&\phantom{00000000}$7$\phantom{00000000}&\phantom{00000000}$8$\phantom{00000000}&\phantom{00000000}$9$\phantom{00000000}\\ 
\midrule
 6 &                    7190 &                       0 &                       0 &                       0\\
 7 &                  469632 &                  144904 &                       0 &                       0\\
 8 &                18521632 &                11990766 &                 3534490 &                       0\\
 9 &               563764626 &               571333104 &               352456980 &               100895667\\
10 &             14578260141 &             20460879142 &             19724988666 &             11802591792\\
11 &            336572753272 &            611002852755 &            814875498464 &            761046274191\\
12 &           7146006383642 &          16057491913346 &          27678691047563 &          35728218734494\\
13 &         142229030736882 &         383591381459308 &         817987730593078 &        1363656535674147\\
14 &        2688967762932621 &        8507419863898968 &       21764920947291140 &       44859920849412082\\
\midrule
\end{tabular}
\caption{Non-orientable genus $g$ maps with $n$ edges}
\label{table:table_last_2}
\end{center}
\end{table}

\begin{table}[h!]
\begin{center}
\footnotesize
\begin{tabular}{c|cccc}
\midrule
$ n\backslash g $ &\phantom{00000000}$10$\phantom{00000000}&\phantom{00000000}$11$\phantom{00000000}&\phantom{00000000}$12$\phantom{00000000}&\phantom{00000000}$13$\phantom{00000000}\\ 
\midrule
10 &              3276228298 &                       0 &                       0 &                       0\\
11 &            441064161280 &            119465644032 &                       0 &                       0\\
12 &          32289907002323 &          18241550095386 &           4827606232542 &                       0\\
13 &        1702134963077638 &        1498266227902654 &         826272319923692 &         214282994249825\\
14 &       72275063245968670 &       87804887825945334 &       75393018976649179 &       40709834082394876\\
\midrule
\end{tabular}
\caption{Non-orientable genus $g$ maps with $n$ edges}
\label{table:table_last_3}
\end{center}
\end{table}

\end{document}